\documentclass[10pt]{article}
\usepackage[includeheadfoot,top=2 cm, bottom=2 cm, left=2 cm, right=2 cm]{geometry}
\usepackage{appendix}
\usepackage{titling}
\usepackage{graphicx}
\usepackage{lipsum}
\usepackage{amsfonts,amsthm}
\usepackage{upgreek}
\usepackage{hyperref}
\usepackage{algorithmicx,algorithm,algcompatible,eqparbox}

\usepackage{amssymb}
\usepackage{multirow}
\usepackage{array}
\usepackage{mathtools}
\usepackage{color}
\usepackage{xcolor}
\usepackage{float}
\usepackage{subcaption}
\usepackage[noabbrev,capitalise]{cleveref}
\hypersetup{
	colorlinks,
	linkcolor={blue!},
	citecolor={blue!},
	urlcolor={blue!}
}

\Crefname{proposition}{Proposition}{Propositions}
\crefname{equation}{}{}
\Crefname{equation}{}{}
\numberwithin{equation}{section}

\usepackage[
sortcites,
backend=bibtex,
hyperref=true,
firstinits=true,
maxbibnames=200,
backref,
backrefstyle=none
]{biblatex} 
\renewbibmacro{in:}{}
\newcommand{\cond}{\textup{cond}}
\newcommand{\fl}{\textup{fl}}

\newcommand{\bDelta}{\mathbf{\Delta}}
\newcommand{\bphi}{\mathbf{\upvarphi}}
\newcommand{\bgamma}{\mathbf{\upgamma}}

\newcommand{\bxi}{\mathbf{\upxi}}
\newcommand{\Frob}{\mathrm{F}}

\newcommand{\ba}{\mathbf{a}}
\newcommand{\bp}{\mathbf{p}}
\newcommand{\bb}{\mathbf{b}}
\newcommand{\bu}{\mathbf{u}}
\newcommand{\be}{\mathbf{e}}
\newcommand{\bv}{\mathbf{v}}
\newcommand{\by}{\mathbf{y}}
\newcommand{\bt}{\mathbf{t}}
\newcommand{\bq}{\mathbf{q}}
\newcommand{\bs}{\mathbf{s}}
\newcommand{\bz}{\mathbf{z}}
\newcommand{\bx}{\mathbf{x}}
\newcommand{\bw}{\mathbf{w}}

\newcommand{\bnull}{\mathbf{0}}
\newcommand{\bTheta}{\mathbf{\Theta}}

\newcommand{\bPhi}{\mathbf{\Phi}}
\newcommand{\bR}{\mathbf{R}}
\newcommand{\bA}{\mathbf{A}}

\newcommand{\bhS}{\widehat{\mathbf{S}}}
\newcommand{\bhR}{\widehat{\mathbf{R}}}
\newcommand{\bhW}{\widehat{\mathbf{W}}}
\newcommand{\bhB}{\widehat{\mathbf{B}}}
\newcommand{\bhQ}{\widehat{\mathbf{Q}}}

\newcommand{\bhX}{\widehat{\mathbf{X}}}
\newcommand{\bhV}{\widehat{\mathbf{V}}}
\newcommand{\bhP}{\widehat{\mathbf{P}}}

\newcommand{\bhH}{\widehat{\mathbf{H}}}
\newcommand{\bhs}{\widehat{\mathbf{s}}}
\newcommand{\bhr}{\widehat{\mathbf{r}}}
\newcommand{\bhx}{\widehat{\mathbf{x}}}
\newcommand{\bhq}{\widehat{\mathbf{q}}}
\newcommand{\bhp}{\widehat{\mathbf{p}}}

\newcommand{\bhw}{\widehat{\mathbf{w}}}

\newcommand{\bX}{\mathbf{X}}

\newcommand{\bT}{\mathbf{T}}

\newcommand{\bH}{\mathbf{H}}

\newcommand{\bW}{\mathbf{W}}

\newcommand{\bS}{\mathbf{S}}
\newcommand{\bY}{\mathbf{Y}}
\newcommand{\bQ}{\mathbf{Q}}
\newcommand{\bC}{\mathbf{C}}
\newcommand{\bP}{\mathbf{P}}
\newcommand{\bI}{\mathbf{I}}

\newcommand{\bV}{\mathbf{V}}
\newcommand{\bF}{\mathbf{F}}
\newcommand{\bPi}{\mathbf{\Pi}}

\newcommand{\bh}{\mathbf{h}}

\newcommand{\bone}{\mathbf{1}}

\newtheorem{lemma}{Lemma}[section]
\newtheorem{proposition}[lemma]{Proposition}
\newtheorem{corollary}[lemma]{Corollary}
\newtheorem{theorem}[lemma]{Theorem}
\newtheorem{remark}[lemma]{Remark}
\newtheorem{definition}[lemma]{Definition}
\newtheorem{assumptions}[lemma]{Assumptions}

\begin{document}

	\title{Randomized Gram-Schmidt process with application to GMRES.}
	\author{Oleg Balabanov\thanks{Alpines,  Inria,  Sorbonne Universit\'e,  Universit\'e de Paris,  CNRS, Laboratoire Jacques-Louis Lions,  F-75012 Paris. Email: oleg.balabanov@inria.fr.}		
		~and~Laura Grigori\thanks{Alpines,  Inria,  Sorbonne Universit\'e,  Universit\'e de Paris,  CNRS, Laboratoire Jacques-Louis Lions,  F-75012 Paris.  Email: laura.grigori@inria.fr.}}
	\date{}
	\maketitle
	\vspace*{-1cm}

\begin{abstract}
	A randomized Gram-Schmidt algorithm is developed for orthonormalization of high-dimensional vectors or QR factorization.
	The proposed process can be less computationally expensive than the classical Gram-Schmidt process while being at least as numerically stable as the modified Gram-Schmidt process. Our approach is based on random sketching, which is a dimension reduction technique consisting in estimation of inner products of high-dimensional vectors by inner products of their small efficiently-computable random {images}, so-called sketches. 
	{In this way, an approximate orthogonality of the full vectors can be obtained by orthogonalization of their sketches.}
	
	The proposed Gram-Schmidt algorithm can provide computational cost reduction in any architecture.   
	{The benefit of random sketching can be amplified  by performing the non-dominant operations in higher precision. In this case the numerical stability can be guaranteed with a working unit roundoff independent of the dimension of the problem.}
	
	The proposed Gram-Schmidt process can be applied to Arnoldi iteration and result in new Krylov subspace methods for solving high-dimensional systems of equations or eigenvalue problems. Among them we chose randomized GMRES method as a practical application of the methodology.    
\end{abstract}


\section{Introduction}
The orthonormalization of a set of high-dimensional vectors serves as basis for many algorithms in numerical linear algebra and other fields of science and engineering. The Gram-Schmidt process (GS) is one of the easiest and most powerful methods to perform this task. 

The numerical stability of the standard implementations of GS, which are the classical Gram-Schmidt algorithm (CGS) and the modified Gram-Schmidt algorithm (MGS), were analyzed in~\parencite{bjorck1967solving,abdelmalek1971round}. The analysis of CGS was improved in~\parencite{giraud2005loss,giraud2005rounding}. In~\parencite{abdelmalek1971round,giraud2005rounding,leon2013gram,smoktunowicz2006note} the authors discussed more sophisticated variants of the GS process and in particular the classical Gram-Schmidt algorithm with re-orthogonalization (CGS2). Versions of GS process well-suited for modern extreme-scale computational architectures were developed in~\parencite{swirydowicz2018low,malard1994efficiency}. 

In this article we propose a probabilistic way to reduce the computational cost of GS process by using the random sketching technique~\parencite{sarlos2006improved,halko2011finding,woodruff2014sketching}. This approach recently became a popular tool for solving high-dimensional problems arising in such fields as theoretical computer science, signal processing, data analysis, model order reduction, and machine learning~\parencite{woodruff2014sketching,vershynin2018high}. The key idea of random sketching technique {relies in the} estimation of inner products of high-dimensional vectors by inner products of their low-dimensional {images through a random matrix}. The random {sketching matrix} is chosen depending on the computational architecture so that it can be efficiently applied to a vector.  In this way, one is able to efficiently embed a set (or a subspace) of high-dimensional vectors, defining the problem of interest, into a low-dimensional space and then {tackle} the problem in this low-dimensional space. In the context of GS process, this implies {orthogonalizing the sketches} rather than high-dimensional vectors. Along with the randomized variant of GS process, here referred to as the randomized Gram-Schmidt process (RGS), we also provide precise conditions on the sketch to guarantee the approximate orthogonality of the output vectors in finite precision arithmetic.  {They rely on the $\varepsilon$-embedding property of the random {sketching matrix} for the subspace spanned by the output vectors. This property is shown to hold {for standard random matrices} with high probability, if the  set of vectors to be orthogonalized is provided a priori.
Furthermore, an efficient procedure for the a posteriori certification of the $\varepsilon$-embedding property is presented. Besides the certification of the output, this procedure can be used for the adaptive selection of the {size} of the random {sketching matrix} or for  improving the robustness of algorithms as depicted in~\cref{rmk:multRGS}.}

Furthermore, we show how the efficiency gains of the RGS algorithm can be amplified by using a multi-precision arithmetic. In particular, it is proposed to perform expensive high-dimensional operations in low precision, {that represents the working precision}, while computing the efficient random projections and low-dimensional operations in high precision. By exploiting statistical properties of rounding errors~\parencite{higham2019new,connolly2020stochastic}, we are able to prove the stability of RGS for the {working} precision unit roundoff independent of the high dimension of the problem.  {Clearly, the presented analysis directly implies stability guarantees also for the unique precision model.}

The randomization {entails}  a possible failure of an algorithm. The probability of this happening, however, is a user-specified parameter that can be chosen very small (e.g., $10^{-10}$) without considerable impact on the overall computational costs. 

One of the uses of the GS process is the computation of an orthonormal basis of a Krylov subspace. This procedure may be used for the solution of high-dimensional eigenvalue problems or systems of equations. In the context of minimal residual methods, such an approach is respectively referred to as the Arnoldi iteration for eigenvalue problems and the generalized minimal residual method (GMRES) for linear systems of equations. For the presentation of these methods, see~\cref{application}. The numerical properties of GMRES were analyzed in~\parencite{drkovsova1995numerical,rozloznik1996numerical,greenbaum1997numerical,paige2006modified}. The usage of variable (or multi) precision arithmetic for Krylov methods, and in particular GMRES, was discussed in~\parencite{van2004inexact,giraud2007convergence,simoncini2007recent,yamazaki2014mixed,carson2021mixed,gratton2019exploiting}.  In the present article we chose the GMRES method  as a practical application of the RGS algorithm. 

The organization of the article is as follows. In~\cref{prem} we describe the basic notations. \Cref{GSprocess} introduces a general GS process and particularizes it to few classical variants. \Cref{RGSA} at first discusses the general idea of the random sketching technique. Then in~\cref{skroundofferrors}, we analyze the rounding errors of a sketched matrix-vector product. A version of GS process, based on random sketching, is proposed in~\cref{RGS}. Its performance in different computational architectures is then studied in~\cref{peranal}. \Cref{GSstability} is devoted to the a priori as well as a posteriori stability analysis of the randomized GS process. \Cref{application} discusses the incorporation of the methodology into the Arnoldi iteration and GMRES algorithms. \Cref{numexp} provides the experimental validation of proposed algorithms.  Finally, \cref{concl} concludes the article. 

{For better presentation most of the proofs of theorems and propositions are provided as supplementary material.}

\subsection{Preliminaries} \label{prem}
Throughout the manuscript we work with real numbers noting that the presented methodology can be naturally extended to complex numbers.

Algebraic vectors are here denoted by bold lowercase letters, e.g, letter $\bx$. For given vectors $\bx_1, \hdots, \bx_k$, we denote matrix $[\bx_1, \hdots, \bx_k]$ by $\bX_k$ (with a bold capital letter) and the $(i,j)$-th entry of $\bX_k$ by $x_{i,j}$ (with a lowercase letter). The notation $\bX_k$ can be further simplified to $\bX$ if $k$ is constant.
Furthermore, we let $[\bX_k]_{(N_1:N_2,M_1:M_2)}$  denote the block of entries $x_{i,j}$ of $\bX_k$ with $(i,j) \in \{N_1, N_1+1, \hdots, N_2 \} \times \{M_1, M_1+1, \hdots, M_2\}$. For a special case of $M_1=M_2$, we denote the vector $[\bX_k]_{(N_1:N_2,M_1:M_2)}$ by simply
$[\bX_k]_{(N_1:N_2,M_1)}$. Moreover, if $\bX_k$ is a vector,  $[\bX_k]_{(N_1:N_2,1)}$ is denoted by $[\bX_k]_{(N_1:N_2)}$.  The minimal and the maximal singular values of $\bX$ are denoted by $\sigma_{min}(\bX)$ and $\sigma_{max}(\bX)$, and the condition number by $\cond(\bX)$. We let $\langle \cdot, \cdot \rangle$ and $\|\cdot\| = \sigma_{max}(\cdot)$ be the $\ell_2$-inner product and $\ell_2$-norm, respectively. $\| \cdot \|_\Frob$ denotes the Frobenius norm.  For two matrices (or vectors) $\bX$ and $\bY$ we say that 
$ \bX \leq \bY, $  
if the entries of $\bX$ satisfy $x_{i,j} \leq y_{i,j}$. Furthermore, for a matrix (or a vector) $\bX$, we denote by $|\bX|$ the matrix $\bY$ with entries $y_{i,j}= |x_{i,j}|$. We also let $\bX^\mathrm{T}$ and $\bX^{\dagger}$ respectively denote the transpose and the Moore--Penrose inverse of $\bX$. Finally, we let $\bI_{k \times k}$ be {the} $k \times k$ identity matrix. 

For a quantity or an arithmetic expression, $X$, we use notation  $\fl(X)$ or $\hat{X}$ to denote the computed value of $X$ with finite precision arithmetic.

\subsection{Gram-Schmidt process} \label{GSprocess}
The Gram-Schmidt process is a method to orthonormalize a set of vectors or compute QR factorization of a matrix. We are concerned with a column-oriented variant of the process.  It proceeds recursively, at each iteration selecting a new vector from the set and orthogonalizing it with respect to the previously selected vectors, as is depicted in~\cref{alg:GS}.  
\begin{algorithm}[h] \caption{Gram-Schmidt process} \label{alg:GS}
	\begin{algorithmic}
		\STATE{\textbf{Given:} $n \times m$ matrix $\bW$, $m \leq n$ }
		\STATE{\textbf{Output}:  $n \times m$ factor $\bQ$ and $m \times m$ upper triangular factor $\bR$.}
		\FOR{ $i = 1:m$} 
		\STATE{1. Compute a projection $\bq_i = \bPi^{(i-1)} \bw_i$  (also yielding $[\bR]_{(1:i-1,i)}$)}.
		\STATE{2. Normalize $\bq_i$ (also yielding $r_{i,i}$)}.
		\ENDFOR
	\end{algorithmic}
\end{algorithm}

The projector $\bPi^{(j)}$ in~\cref{alg:GS}  is usually taken as approximation to the ($\ell_2$-)orthogonal projector $\bI_{n \times n} - \bQ_j (\bQ_j)^\dagger$ onto $\mathrm{span}(\bQ_{j})^{\perp}$, {$1 \leq j \leq m-1$}. If~{\Cref{alg:GS}} is used with infinite precision arithmetic, then considering
\begin{equation} \label{eq:classproj}
\bPi^{(j)} = \bI_{n \times n} - \bQ_j (\bQ_j)^\mathrm{T}
\end{equation}
will produce an exact QR factorization of $\bW$. This fact can be shown by induction. In short, we can show that $\bPi^{(i-1)}$ being an orthogonal projector implies $\bQ_{i}$ being an orthonormal matrix, which in its turn implies that $\bPi^{(i)}$ is an orthonormal projector. \Cref{alg:GS} with the choice~\cref{eq:classproj} is referred to as the classical Gram-Schmidt process. With finite precision arithmetic, however, matrix $\bQ_{j}$ can be guaranteed to be orthonormal only approximately. {This} can make the classical GS algorithm suffer from numerical instabilities.  In particular, in this case the orthogonality of Q factor, measured by $\|\bI_{m \times m} - \bQ^\mathrm{T} \bQ \|$,  can grow as $\cond(\bW)^2$ {or more, depending on the method used for normalization~\parencite{bjorck1967solving,smoktunowicz2006note}.}

Besides the projector~\cref{eq:classproj}, there are a couple of other standard choices for $\bPi^{(j)}$. They can yield a better numerical stability but require more computational cost {in terms of flops, storage consumption, scalability or amount of communication between processors}. Modified Gram-Schmidt algorithm {uses the} projector
$$\bPi^{(j)} = (\bI_{n \times n} - \bq_j (\bq_j)^\mathrm{T}) (\bI_{n \times n} - \bq_{j-1} (\bq_{j-1})^\mathrm{T}) \hdots (\bI_{n \times n} - \bq_1 (\bq_1)^\mathrm{T}),~1\leq j \leq m-1.$$
In this case the orthogonality measure of $\bQ$ depends only linearly on $\cond(\bW)$~\parencite{bjorck1967solving}. Numerical stability of the modified GS algorithm is sufficient for most applications and is often considered as benchmark for characterizing {the} stability of algorithms for {orthogonalizing a set of vectors or computing a QR factorization.} Another choice for $\bPi^{(j)}$ is 
$$\bPi^{(j)} = (\bI_{n \times n} - \bQ_j (\bQ_j)^\mathrm{T})(\bI_{n \times n} - \bQ_j (\bQ_j)^\mathrm{T}),~1\leq j \leq m-1,$$
which results in a so-called classical Gram-Schmidt process with re-orthogonalization (CGS2). This projector can be shown to yield a similar (or better) stability as the modified GS process.

In this work we are concerned with a scenario when $\bW$ is a large matrix with a moderate number of columns, i.e., when $m \ll n$.   For this situation, we propose a new {\emph{randomized}} projector $\bPi^{(j)}$ that can yield more efficiency than {CGS} process, while providing no less numerical stability than {MGS} process. {Unlike standard approaches, our RGS algorithm  provides a Q factor that is not $\ell_2$-orthogonal even under exact arithmetic, but that is very well-conditioned with very high probability. 
This property is sufficient for a number of applications.
For instance, as is shown in~\cref{RGMRES}, a small condition number of the Q factor guarantees an almost optimal convergence of the GMRES solution. For other cases, the Q factor produced by RGS algorithm should be post processed with a Cholesky QR.} 

\section{Randomized Gram-Schmidt algorithm} \label{RGSA}
\subsection{Introduction to random sketching} \label{skintro}

Let $\bTheta \in \mathbb{R}^{k \times n}$, with $k \ll n$, be a sketching matrix. This matrix shall be seen as an embedding of subspaces of $\mathbb{R}^{n}$ into subspaces of $\mathbb{R}^{k}$ and is therefore referred to as a $\ell_2$-subspace embedding. The $\ell_2$-inner products between vectors in subspaces of $\mathbb{R}^{n}$ are estimated by $$ \langle \cdot, \cdot \rangle \approx \langle \bTheta \cdot, \bTheta \cdot \rangle. $$
For a given (low-dimensional) subspace of interest $V \subset \mathbb{R}^{n}$, the quality of such an estimation can be characterized by the following property of $\bTheta$.
\begin{definition}
	For $\varepsilon<1$, the sketching matrix $\bTheta \in \mathbb{R}^{k \times n}$ is said to be {an} $\varepsilon$-subspace embedding for $V \subset \mathbb{R}^{n}$, if we have
	\begin{equation} \label{eq:isometry}
	\forall \bx,\by \in V,~~ | \langle \bx, \by \rangle - \langle \bTheta \bx, \bTheta \by \rangle | \leq \varepsilon \| \bx \| \|\by\|.  
	\end{equation} 
\end{definition}

Let $\bV$ be a matrix whose columns form a basis for $V$. To ease presentation in next sections, {an} $\varepsilon$-subspace embedding for $V$ shall be often referred to simply as {an} $\varepsilon$-embedding for $\bV$.
\begin{corollary} \label{thm:epscond}
	If $\bTheta \in \mathbb{R}^{k \times n}$ is {an} $\varepsilon$-embedding for $\bV$, then the singular values of $\bV$ are bounded by $$ (1+\varepsilon)^{-1/2} \sigma_{min}(\bTheta\bV) \ \leq \sigma_{min}(\bV) \leq \sigma_{max}(\bV) \leq  (1-\varepsilon)^{-1/2} \sigma_{max}(\bTheta\bV).$$
	\vspace*{-0.5cm}	
	\begin{proof}
		Let $\ba \in \mathbb{R}^{\mathrm{dim}(V)}$ be an arbitrary vector and $\bx = \bV \ba$.
		By definition of $\bTheta$, 
		$$ (1+\varepsilon)^{-1} \| \bTheta \bx \|^2  \leq \| \bx \|^2 \leq (1-\varepsilon)^{-1} \| \bTheta \bx \|^2, \text{ which implies that }$$ 
		$$ (1+\varepsilon)^{-1/2} \| \bTheta \bV \ba \|  \leq \| \bV \ba \| \leq (1-\varepsilon)^{-1/2} \| \bTheta \bV \ba \|.$$	
		The statement of proposition then follows by using definitions of the minimal and the maximal singular values of a matrix.		
	\end{proof}
\end{corollary}
\Cref{thm:epscond} implies that to make the condition number of matrix $\bV$ close to $1$, it can be sufficient to orthonormalize small sketched matrix $\bTheta \bV$. This observation serves as basis for the randomized GS process in~\cref{RGS}.  Note that the orthogonalization of $\bTheta \bV$ with respect to the $\ell_2$-inner product is equivalent to orthonormalization of $\bV$ with respect to {the}  product $\langle \bTheta \cdot , \bTheta \cdot \rangle$. Note also that in our applications there will be no practical benefit of considering very small values for $\varepsilon$. The usage of $\varepsilon \leq 1/2$ or $\varepsilon \leq 1/4$ will be sufficient.

We {here} proceed with sketching matrices that do not require any a priori knowledge of $V$ to guarantee~\cref{eq:isometry}. Instead, $\bTheta$ {is} generated from a carefully chosen distribution such that it satisfies~\cref{eq:isometry} for any low-dimensional subspace with high probability.
\begin{definition} \label{def:oblemb}
	The sketching matrix $\bTheta \in \mathbb{R}^{k \times n}$ is called a $(\varepsilon,\delta,d)$ oblivious $\ell_2$-subspace embedding, if it is {an} $\varepsilon$-embedding for any fixed $d$-dimensional subspace $V \subset \mathbb{R}^{n}$ with probability at least $1-\delta$.
\end{definition}
In general, {such \textit{oblivious subspace embeddings}} with high probability have a bounded norm, as is shown in~\cref{prop:Cn2sqrtn}.
\begin{corollary} \label{prop:Cn2sqrtn}
	If $\bTheta \in \mathbb{R}^{k \times n}$ is a $(\varepsilon, \delta/n, 1)$ oblivious $\ell_2$-subspace embedding, then with probability at least $1-\delta$, we have
	$$\|\bTheta\|_\Frob \leq\sqrt{(1+\varepsilon)n}.$$
	\vspace*{-0.5cm}	
	\begin{proof}
		It directly follows from~\cref{def:oblemb} and the union bound argument that $\bTheta$ is an $\varepsilon$-embedding for each canonical (Euclidean) basis vector. This implies that the $\ell_2$-norms of the columns of $\bTheta$ are bounded from above by $\sqrt{1+\varepsilon}$. The statement of the corollary then follows immediately. 
	\end{proof}
\end{corollary}

There are several distributions that are known to satisfy the $(\varepsilon,\delta,d)$ oblivious $\ell_2$-subspace embedding property when $k$ is sufficiently large. The standard examples include  Gaussian, Rademacher distributions, sub-sampled randomized Hadamard and Fourier transforms, CountSketch matrix and more~\parencite{halko2011finding,woodruff2014sketching,achlioptas2003database}. In this work we shall rely on Rademacher matrices and partial SRHT (P-SRHT). A (rescaled) Rademacher matrix has i.i.d. entries equal to $\pm1/\sqrt{k}$ with probabilities $1/2$.  The efficiency of multiplication by Rademacher matrices
can be attained due to proper exploitation of computational architectures.  For instance the products of Rademacher matrices with vectors can be implemented with standard SQL primitives and are embarrassingly parallelizable. For $n$ being a power of $2$, SRHT is defined as a product of diagonal matrix of random signs with Walsh-Hadamard matrix, followed by  uniform sub-sampling matrix and scaling factor $1/\sqrt{k}$. Random sketching with SRHT can improve efficiency in terms of number of flops. Products of SRHT matrices with vectors require only $n \log_2(n)$ flops  using  the  fast  Walsh-Hadamard  transform or $2n\log_2(k+ 1)$ flops using the procedure in~\parencite{ailon2009fast}.  P-SRHT is used instead of SRHT when $n$ is not a power of $2$ and is defined as the first $n$ columns of SRHT matrix of size $s$, were $s$ is the power of $2$ such that $ n \leq s < 2n$. Furthermore, for both (P-)SRHT and Rademacher matrices a seeded random number generator can be utilized to allow efficient {storage and application} of $\bTheta$.  This is particularly important for limited-memory and distributed computational architectures. It follows that the rescaled Rademacher distribution with 
\begin{subequations} \label{eq:ell_2bounds}
	\begin{align}
	k &\geq 7.87 \varepsilon^{-2}( {6.9} d + \log (1/\delta)), \\
	\intertext{and the P-SRHT distribution with} 
	k &\geq 2( \varepsilon^{2} - \varepsilon^3/3)^{-1} \left (\sqrt{d}+ \sqrt{8 \log(6 n/\delta)} \right )^2 \log (3 d/\delta),
	\end{align}
\end{subequations}
respectively, are $(\varepsilon, \delta, d)$ oblivious $\ell_2$-subspace embeddings~\parencite{balabanov2019randomized}. We see that the bounds~\cref{eq:ell_2bounds} are independent or only logarithmically dependent on the dimension $n$ and probability of failure, and are proportional to the low dimension $d$. {This implies that one can use $\bTheta$ of a small size even for very large problems and very small probabilities of failure.} 

\subsection{Rounding errors in a sketched matrix-vector product} \label{skroundofferrors}

Let us fix a realization of an oblivious $\ell_2$-subspace embedding  $\bTheta \in \mathbb{R}^{k \times n}$ of sufficiently large size, and consider a matrix-vector product
$$ \bx = \bY \bz,\text{ with $\bY \in \mathbb{R}^{n\times m},~\bz \in \mathbb{R}^{m} $} , $$
computed in finite precision arithmetic with unit roundoff $u <0.01/m$. Note that elementary linear algebra operations on vectors such as addition or multiplication by constant can be also viewed as matrix-vector products.  Define rounding error vector $\bDelta \bx  = \bhx - \bx$.  The standard worse-case scenario rounding analysis provides an upper bound for $\bDelta \bx$ of the following form~\parencite{higham2002accuracy}:
\begin{equation} \label{eq:roundoffYz}
|\bDelta \bx| \leq  \bu
\end{equation}
In a general case, the vector $\bu$ can be taken as\footnote{We here used the fact that $\frac{mu}{1-m u} \leq 1.02 mu$.} 
\begin{equation} 
\bu = 1.02 m u |\bY| |\bz|.
\end{equation}
In some situations, e.g., if the matrix $\bY$ is sparse, this bound can be improved. Here, we are particularly interested in the case, when $\bY = a \bI_{n \times n}$, i.e., when $\bY \bz$ represents a multiplication of $\bz$ by a constant, and when $\bY \bz = \bY' \bz'+\bh$, i.e., when it represents a sum of a matrix-vector product with a vector. Then in the first case, one can take $$ \bu = u |a \bz|,$$ and in the second case,\footnote{We have by the standard worse-case scenario analysis, $$ |\bDelta \bx|  \leq \frac{(m-1)u}{1-(m-1)u} |\bY'||\bz'| + u \left ((1+\frac{(m-1)u}{1-(m-1)u}) |\bY'||\bz'|+ |\bh| \right ) \leq 1.02 u (|\bh| + m |\bY'| |\bz'|)$$} $$\bu =  1.02 u (|\bh| + m |\bY'| |\bz'|).$$

Let us now address bounding the rounding error of the sketch $\bTheta \bhx$. This will become particularly handy in~\cref{GSstability} to simplify stability analysis of the proposed in~\cref{RGS} randomized GS algorithm. We {here} seek a bound of the form
\begin{equation} \label{eq:numepsembedding}
\| \bTheta \bDelta \bx\| \leq D \| \bu \|, 
\end{equation}
where $D$ is a coefficient possibly depending on $n$. Clearly, we have
\begin{equation} \label{eq:numepsembedding2}
\| \bTheta \bDelta \bx\| \leq \|\bTheta\| \| \bu \|, 
\end{equation}
{which combined with~\cref{prop:Cn2sqrtn} implies that with high probability the relation~\cref{eq:numepsembedding} holds with $D = \mathcal{O}(\sqrt{n})$}.\footnote{If $\bTheta$ is a $(\varepsilon, \delta/n, 1)$ oblivious subspace embedding, then we have $D = \sqrt{1+\varepsilon}\sqrt{n}$ with probability at least $1-\delta$.} 

{Next we notice that taking $D = \mathcal{O}(\sqrt{n})$ accounts for a very improbable worse-case scenario and is pessimistic in practice. If $\bx$ is independent of $\bTheta$, then with high probability the relation~\cref{eq:numepsembedding} holds for $D = \mathcal{O}(1)$.}\footnote{If $\bTheta$ is a $(\varepsilon, \delta, 1)$ oblivious subspace embedding, then we have $D = \sqrt{1+\varepsilon}$ with probability at least $1-\delta$.} 
Furthermore, (as it is argued below in details) one can expect the relation~\cref{eq:numepsembedding}, with $D = \mathcal{O}(1)$, to hold for  $\bTheta$ of moderate size even when $\bx$ depends on $\bTheta$ (i.e., when $\bY$ and $\bz$ are chosen depending on $\bTheta$), since the rounding error vector $\bDelta \bx$ should in practice have only a minor correlation with $\bTheta$. {This property can be viewed as a sketched version of the standard ``rule of thumb'' stating that in practice one can reduce the worse-case scenario error constants (e.g., constant $\gamma_n = \frac{nu}{1-nu}$ in~\parencite{higham2002accuracy}) by a factor of $\sqrt{n}$.} It has important meaning in the context of (oblivious) randomized algorithms: the sketching step does not in practice multiply the rounding errors by a factor depending on $n$. In other words, with random sketching one is able to efficiently reduce the dimension of the problem without a loss of numerical precision. 

To provide a precise guarantee that~\cref{eq:numepsembedding} holds for $D = \mathcal{O}(1)$ we shall need to explore the properties of $\bDelta \bx$ as a vector of rounding errors. For this we shall consider a {probabilistic} rounding model, where 
\begin{itemize}
	\item {the rounding errors $\xi$ due to each elementary arithmetic operation $x~\mathrm{op}~y$, i.e., 
		$$\xi = \frac {\fl(x~\mathrm{op}~y) - (x~\mathrm{op}~y)} {(x~\mathrm{op}~y)},~\text{with } \mathrm{op} = +,-,*,/,$$	 
		are bounded random variables possibly depending on each other, but are independently centered (i.e., have zero mean).}	
	\item the computation of each entry of $\bhx$ is done independently of other entries. In other words, the entries of $\bDelta \bx$  are drawn independently of each other.  
\end{itemize}
{This model corresponds to~\parencite[Model 4.7]{connolly2020stochastic}. Its particular case is the so-called
stochastic rounding model (see~\parencite{connolly2020stochastic}), which recently gained  attention in {the} machine learning community to improve the accuracy and the efficiency of training neural networks.} The analysis of standard numerical linear algebra algorithms, and in particular, the rigorous foundation of the ``rule of thumb'', with the {probabilistic} rounding model is provided in~\parencite{higham2019new, higham2020sharper, connolly2020stochastic,ipsen2019probabilistic}. Note that the used here rounding model does not assume the rounding errors to be independent random variables as in~\parencite{higham2019new,higham2020sharper,ipsen2019probabilistic}, but only mean-independent with zero mean, which is a weaker and more realistic assumption, as {it} is argued in~\parencite{connolly2020stochastic}.

Let us deduce that under the described probabilistic model, the vector $\bDelta \bx$ has entries that are independent centered random variables. It then follows from~\cref{thm:thetadeltax} that $\bTheta$ shall satisfy~\cref{eq:numepsembedding} with $D = \mathcal{O}(1)$  with probability at least $1-2\delta$, if $\bTheta$ is a  $(\varepsilon,\binom{n}{d}^{-1}\delta, d)$, with $d = \mathcal{O}(\log (1/\delta))$, oblivious $\ell_2$-subspace embedding. According to~\cref{eq:ell_2bounds}, this property is satisfied if $\bTheta$ is a Rademacher matrix with {$\mathcal{O}(\log(n)\log(1/\delta))$} rows or P-SRHT matrix with $\mathcal{O}(\log^2(n)\log^2(1/\delta))$ rows.

\begin{theorem} \label{thm:thetadeltax}
	Draw a realization $\bTheta \in \mathbb{R}^{k \times n}$ of  $(\varepsilon/4,\binom{n}{d}^{-1} \delta, d)$ oblivious $\ell_2$-subspace embedding, with $d = 4.2 c^{-1} \log (4/\delta)$, where $c \leq 1$ is some universal constant.
	Let $\bphi \in \mathbb{R}^{n}$ be a vector with entries that are independent random variables from  distributions that can depend on $\bTheta$. If $\bphi$ has zero mean, i.e., {$E(\bphi|\bTheta) = \bnull$}, and	$|\bphi| \leq \bgamma$ for some vector $\bgamma \in \mathbb{R}^{n}$, then 
	\begin{equation} \label{eq:thetadeltax11}
	|\| \bphi\|^2 - \| \bTheta \bphi\|^2| \leq \varepsilon \| \bgamma \|^2,
	\end{equation}	
	holds with probability at least $1-2\delta$. 
	\begin{proof} 
		The proof of~\cref{thm:thetadeltax} will rely on the following property of $\bTheta$: 
		\begin{equation} \label{eq:rip}
		(1-\varepsilon)\|\ba\|^2 \leq \|\bTheta \ba\|^2 \leq (1+\varepsilon)\|\ba\|^2 , \text{ for all $d$-sparse vectors $\ba \in \mathbb{R}^n$},
		\end{equation}
		called Restricted Isometry Property of level $\varepsilon$ and order $d$, or simply $(\varepsilon, d)$-RIP.\footnote{A vector $\ba \in \mathbb{R}^n$ is called $d$-sparse if it has at most $d$ nonzero entries.} This is a well-known fact that oblivious $\ell_2$-subspace embeddings satisfy the RIP property with high probability, as is shown in~\cref{thm:obl2rip}. The statement of~\cref{thm:thetadeltax} then follows by combining~\cref{thm:obl2rip} with~\cref{thm:thetadeltax0}, and the union bound argument. \qed
	\end{proof}
\end{theorem}

\begin{corollary} \label{thm:thetadeltax2}
	Consider the {probabilistic} rounding model. If $\bTheta$ is a $(\varepsilon/4,\binom{n}{d}^{-1}\delta, d)$ oblivious $\ell_2$-subspace embedding, with $d = 4.2 c^{-1} \log (4/\delta)$, where $c \leq 1$ is some universal constant, then the bound~\cref{eq:numepsembedding} holds with $D = \sqrt{1+\varepsilon}$ with probability at least $1-2\delta$. 
\end{corollary}

\begin{remark}
	The universal constant $c$ in~\cref{thm:thetadeltax,thm:thetadeltax2} is same as that in the Hanson-Wright inequality (see~\parencite[Theorem 6.2.1]{vershynin2018high}). It can be shown that this constant is greater than $1/64$~\parencite{pollard17}.
\end{remark}

\begin{proposition} \label{thm:obl2rip}
	$(\varepsilon,\binom{n}{d}^{-1}\delta, d)$ oblivious $\ell_2$-subspace embedding $\bTheta \in \mathbb{R}^{k \times n}$ satisfies $(\varepsilon, d)$-RIP with probability at least $1-\delta$.
	\begin{proof}
		Let $\mathcal{B}$ denote the canonical (Euclidean) basis for $\mathbb{R}^n$. It follows directly from the definition of $\bTheta$ and the union bound argument that $\bTheta$ is an $\varepsilon$-embedding for all subspaces spanned by $d$ vectors from $\mathcal{B}$, simultaneously, with probability at least $1-\delta$. Since, every $d$-sparse vector $\ba \in \mathbb{R}^n$ belongs to a subspace spanned by $d$ vectors from $\mathcal{B}$, we conclude that $\bTheta$ satisfies~\cref{eq:rip} with probability at least $1-\delta$.
	\end{proof}
\end{proposition}
\begin{theorem} \label{thm:thetadeltax0}
	Let $\bTheta \in \mathbb{R}^{k \times n}$ be a matrix satisfying $(\varepsilon/4, 2d)$-RIP with $d = 2.1 c^{-1} \log (4/\delta)$, where $c \leq 1$ is some universal constant. Let $\bphi \in \mathbb{R}^{n}$ be a vector with entries that are independent random variables. If $\bphi$ has zero mean and	$|\bphi| \leq \bgamma$ for some vector $\bgamma \in \mathbb{R}^{n}$, then 
	$$ |\| \bphi\|^2 - \| \bTheta \bphi\|^2| \leq \varepsilon \| \bgamma \|^2, $$	
	holds with probability at least $1-\delta$. 
	\begin{proof}
		See supplementary material.
	\end{proof}
\end{theorem}

\subsection{Randomized Gram-Schmidt process} \label{RGS}

Consider the variants of \cref{alg:GS} where the projector $\bPi^{(i-1)}$ has the form 
\begin{equation} \label{eq:GSvariants}
\bPi^{(i-1)} \bw_{i} = \bw_{i} -\bQ_{i-1} \bx,
\end{equation}
with $\bx  = \bR_{(1:i-1,i)}$ computed from $\bQ_{i-1}$ and $\bw_{i}$. 
The classical methods proceed with taking $\bx$ as an approximation of $\bQ_{i-1}^\dagger \bw_{i}$, or equivalently, as an approximate solution to the following least-squares problem:
\begin{equation} \label{eq:GSminimization}
\min_{\by} \|\bQ_{i-1} \by - \bw_{i} \|. 
\end{equation}
The stability of~\cref{alg:GS} in this case can be directly linked to the accuracy of $\bx$.
The CGS and MGS algorithms belong to the aforementioned category of GS processes with $\bx$ taken as, respectively, $\bx = \bQ_{i-1}^\mathrm{T} \bw_{i}$ and $\bx = \bT_{i-1} (\bQ_{i-1}^\mathrm{T} \bw_{i})$, for some triangular matrix $\bT_{i-1}$\parencite{swirydowicz2018low}. The connection of CGS and MGS with solving~\cref{eq:GSminimization} was explored in~\parencite{ruhe1983numerical}. A similar formulation was also used in~\parencite{barlow2005improved}. 


In this work we develop new variants of GS process that satisfy~\cref{eq:GSvariants}, but this time produce the output Q factor orthonormal with respect to the sketched product $\langle \bTheta \cdot , \bTheta \cdot \rangle$ rather than $\ell_2$-inner product as in standard methods. {Thus, the Q factor is no longer $\ell_2$-orthonormal even in exact arithmetic, though, according to~\cref{thm:epscond}, it has a small ($\ell_2$-)condition number and yields a reduced computational cost.} 
{Such factorization} corresponds to taking $\bx$ in~\cref{eq:GSvariants} as an approximation of  $(\bTheta \bQ_{i-1})^\dagger (\bTheta \bw_{i})$, or equivalently, a minimizer of the sketched residual:
$$\min_{\by} \|(\bTheta\bQ_{i-1}) \by - \bTheta \bw_{i} \|.$$  
{Furthermore, the normalization of $\bq_i$ at Step 2 of~\cref{alg:GS} has to be performed accordingly: $\bq_i= \bq_i/\|\bTheta \bq_i \|$.}  
We see that unlike in standard methods, here the computation of $\bx$ requires only (efficient) evaluation of random projections and operations on small vectors and matrices with no standard operations on high-dimensional vectors.  The GS process with such a projector is depicted in~\cref{alg:RGM}. 

{In general, stability of \cref{alg:RGM} directly depends on the  
accuracy and stability of the least-squares solver used in Step 2. One should prioritize least-squares solvers that are as accurate and stable as possible. They  can be based on Givens rotations or Householder transformation as is considered in our stability analysis (see~\cref{stability}).  Such standard solvers should yield a negligible computational cost when matrix $\bS_{m} = \bTheta\bQ_{m}$ is sufficiently small, which happens in most applications. 
However, when $\bS_{m}$ is of moderate size, the least-squares solution with standard methods can entail a considerable computational cost and has to be avoided.
In such cases, by using the fact that $\bS_{i-1}$ is approximately orthonormal, one can compute $\bx = [\bR]_{(1:i-1,i)}$ from the normal equation:
$$ (\bS_{i-1})^\mathrm{T}  \bS_{i-1} \bx = (\bS_{i-1})^\mathrm{T} \bp_i, $$
with several Richardson iterations $\bx \leftarrow  \bx + \bS_{i-1}^\mathrm{T} (\bp_i-\bS_{i-1} \bx)$ requiring a minor computational cost. The resulting algorithm can be viewed as a sketched version of the classical GS process with re-orthogonalizations. The case $\bx = (\bS_{i-1})^\mathrm{T} \bp_i$  with only one Richardson iteration exactly corresponds to the orthogonalization of $\bW$ with respect to $\langle \bTheta \cdot, \bTheta \cdot \rangle$ with the classical Gram-Schmidt process defined for a general inner product.  Moreover, instead of using the Richardson iterations we could also compute $\bx$  by orthogonalizing $\bp_i$  to $\bS_{i-1}$ with an MGS step.	In this case~\cref{alg:RGM} would correspond to orthogonalization of $\bW$ with respect to $\langle \bTheta \cdot , \bTheta \cdot \rangle$ with the modified GS process.
The ways for efficient and stable solution of the sketched least-squares problem are addressed in details in our subsequent work on the block variants of the randomized GS process.  }

\begin{algorithm}[h] \caption{Randomized Gram-Schmidt algorithm (RGS)} \label{alg:RGM}
	\begin{algorithmic}
		\STATE{\textbf{Given:} $n \times m$ matrix $\bW$, and $k \times n$ matrix $\bTheta$, $m \leq k \ll n$.} 
		\STATE{\textbf{Output}:  $n \times m$ factor $\bQ$ and $m \times m$ upper triangular factor $\bR$.}
		\FOR{$i = 1:m$} 
		\STATE{1. Sketch $\bw_i$: $\bp_i = \bTheta \bw_i$.} \COMMENT{macheps: $u_{fine}$}
		\STATE{2. Solve $k \times (i-1)$ least-squares problem: \\~~~~~~~~~~~~~~~~~~~~~~~~~~~{$[\bR]_{(1:i-1,i)} =  \arg \min_\by \| \bS_{i-1} \by - \bp_i \|.$} }
		\COMMENT{macheps: $u_{fine}$}
		\STATE{3. Compute projection of $\bw_i$:  $\bq'_i = \bw_i - \bQ_{i-1}[\bR]_{(1:i-1,i)}$.} \COMMENT{macheps: $u_{crs}$}
		\STATE{4. Sketch $\bq_i'$: $\bs'_i = \bTheta \bq_i'$.}  \COMMENT{macheps: $u_{fine}$}
		\STATE{5. Compute the sketched norm $r_{i,i} = \|\bs'_i\|$.} \COMMENT{macheps: $u_{fine}$}
		\STATE{6. Scale vector $\bs_i = \bs'_i/r_{i,i}$.} \COMMENT{macheps: $u_{fine}$}
		\STATE{7. Scale vector $\bq_i = \bq'_i/r_{i,i}$.} \COMMENT{macheps: $u_{fine}$}
		\ENDFOR
		\STATE{8. (Optional) compute $\Delta_m  = \|\bI_{m \times m} - \bS_m^\mathrm{T} \bS_m \|_\Frob$ and $\tilde{\Delta}_m = \frac{\|\bP_m - \bS_m \bR_m \|_\Frob}{\|\bP_m\|_\Frob}$.\\~~~~Use~\cref{thm:maintheorem1} to certify the output.} \COMMENT{macheps: $u_{fine}$}	
	\end{algorithmic}
\end{algorithm}
{At $i=1$ of~\cref{alg:RGM} we used a conventional notation that $[\bR]_{(1:i-1,i)}$ is a $0$-by-$1$ matrix and $\bQ_{i-1}$ is $n$-by-$0$ matrix, implying that $\bq_i' = \bw_i$ and $\bs_i' = \bp_i$.} \Cref{alg:RGM} is executed with  a multi-precision finite arithmetic  with two unit roundoffs: a coarse one $u_{crs}$, and a fine one $u_{fine}$, $u_{fine} \leq u_{crs} \leq 0.01/m$. The roundoff $u_{crs}$ {represents the \textit{working precision}} and is used for standard operations on high-dimensional vectors {and matrices} in Step 3, which is the most expensive computation in the algorithm.  {This precision is also used for storage of large matrices $\bQ$ and $\bW$.}   All other {(inexpensive)} operations in~\cref{alg:RGM} are performed {and accumulated} with a fine roundoff $u_{fine}$.   We chose a multi-precision model rather than a unique precision one to show an interesting property of the algorithm:  that one may guarantee stability of~\cref{alg:RGM} by performing standard operations on high-dimensional vectors (i.e., Step 3) with unit roundoff $u_{crs}$ independent of $n$ (and $k$). This feature of RGS process can have a particular importance for extreme-scale problems. {Clearly, the results from this paper can be also used for the analysis of~\cref{alg:RGM} executed with {unique} unit roundoff $u_{fine}$. The stability guarantees in such a case can be derived from~\cref{GSstability,application} simply by introducing a fictitious unit roundoff $u_{crs} = F(m,n) u_{fine}$, where $F(m,n)$ is a low-degree polynomial, and looking at~\cref{alg:RGM} as that it is executed with multi-precision arithmetic with unit roundoffs $u_{fine}$ and $u_{crs}$.}
\begin{remark}
	{In Step 4 of~\cref{alg:RGM}, the sketch of $\bq_i'$ could be computed also as $\bs_i' = \bp_i - \bS_{i-1} [\bR]_{(1:i-1,i)}$ {instead of} $\bs_i' = \bTheta \bq_i'$. Our experiments, however, revealed that this way is less stable.}
\end{remark}

\subsection{Performance analysis} \label{peranal}
Let us now characterize the efficiency of~\cref{alg:RGM} executed {in} different computational architectures. The performance analysis  is done through comparison to CGS. The CGS algorithm is the most efficient from the standard (column-oriented) algorithms for orthogonalization of a set of vectors. 
{It requires nearly half as many flops and  {synchronizations between processors} than CGS2 and, unlike MGS, it can be implemented using matrix-vector operations, i.e., level-2 Basic Linear Algebra Subprograms (BLAS).}

By assuming that $m \leq k \ll n$, we shall neglect the cost of operations on sketched vectors and matrices in~\cref{alg:RGM}. Then the computational cost is characterized by evaluation of random sketches at Steps 1 and 4, and  computation of the projection of $\bw_i$ at Step 3. Moreover, at Steps 1 and 4 we let the sketching matrix $\bTheta$ to be chosen depending on each particular situation to yield the most efficiency. 

The RGS algorithm can be beneficial in terms of the classical metric of efficiency, which is the number of flops. If $\bTheta$ is taken as SRHT, then the random projections at Steps 1 and 4, require in total (no more than) $4 n \log(k+1)$ flops at each iteration. For sufficiently large $i$, this cost is much less than the cost of Step $3$ that is nearly  $2ni$ flops. The CGS requires more than $4ni$ flops at each iteration and therefore it is nearly as twice more expensive as RGS.  Furthermore, the flops at Step $3$ of the RGS algorithm can be done in low precision, which can make~\cref{alg:RGM} even more efficient. 

Both CGS and RGS algorithms can be implemented by using BLAS-2 routines for high-dimensional operations. The CGS algorithm in such implementation, however, entails (at least) two passes over the basis matrix $\bQ_{i-1}$ at iteration $i$, $2 \leq i \leq m$. \Cref{alg:RGM}, on the other hand, at each iteration requires only one pass over~$\bQ_{i-1}$ and two applications of $\bTheta$. The applications of $\bTheta$ can be performed by utilizing a seeded random number generator with negligible storage costs. Consequently, in this case RGS can be more pass-efficient than CGS. Furthermore, the matrix $\bQ_{i-1}$ in the RGS algorithm can be maintained in lower precision, {and still yield similar (or better) accuracy} than the CGS algorithm, which can amplify the storage reduction.

To characterize the performance of RGS in parallel/distributed computational architecture, we consider the situation when the columns of $\bW$ are provided recursively as 
$$ \bw_{i+1} = \bA \bq_{i},~1 \leq i \leq m-1, $$
where $\bA$ is a $n \times n$ matrix. 
This, for instance, happens in Arnoldi algorithm for computing an orthonormal basis of a Krylov subspace, which is the core ingredient of GMRES algorithm (see~\cref{application} for details). We here assume that the high-dimensional matrix $\bA$ and the vectors $\bq_{i}$ are distributed among processors using block row-wise partitioning (possibly with overlaps). {This is a standard situation when $\bA$ is obtained from discretization of a PDE.} It is then assumed that the computation of the  matrix-vector product $ \bw_{i+1} = \bA \bq_i$ requires communication only between neighboring processors and has a minor impact on the overall communication cost. We also assume that along with the local matrices and vectors on each processor are also maintained copies of the sketches $\bS_i$, $\bs_i'$ and $\bp_i$, and matrix $\bR_i$.

Next we notice that the utilization of a seeded random number generator can allow efficient access to any block of $\bTheta$ with a minor computational cost and, in particular, with absolutely no communication.  The computation of the sketch $\bp_i = \bTheta \bA \bq_{i-1}$ in Step 1 then requires only one global synchronization. The computation of the sketch $\bTheta \bq'_i$ in Step 4 of~\cref{alg:RGM}  requires an additional synchronization, which implies in total two global synchronizations at each iteration of~\cref{alg:RGM}. This communication cost is the same as of the classical implementation of CGS. {In~\parencite[Section 4]{kim1992efficient}} is described a way to reduce the communication cost of CGS algorithm to only one synchronization per iteration. This technique may also be applied to the RGS algorithm by  
incorporating a lag into Steps 5-7 of~\cref{alg:RGM}. More specifically, at iteration $i$, we can compute two sketches $$ \bs'_i = \bTheta \bq'_i \textup{ and } \bp'_{i+1} = \bTheta \bA \bq'_{i},~1 \leq i \leq m-1,$$ simultaneously,  by utilizing only one global synchronization, and then perform the normalizations: $r_{i,i} = \| \bs'_i \|$, $\bs_i = \bs'_i/r_{i,i}$, $\bq_i = \bq'_i/r_{i,i}$ and $\bp_{i+1} = \bp'_{i+1}/r_{i,i}$. The communication cost of such implementation of the RGS algorithm then becomes only one global synchronization per iteration. We conclude that RGS and CGS should have similar {number of required synchronizations} in parallel/distributed computational architecture.

\section{Stability of randomized Gram-Schmidt process} \label{GSstability}

In this section we provide stability analysis of~\cref{alg:RGM}. It is  based on the following assumptions that hold with high probability if $\bTheta$ is an oblivious subspace embedding of sufficiently large size. 

First, we assume that $\bTheta$ satisfies~\cref{eq:Thetanorm}. According to~\cref{prop:Cn2sqrtn}, this property holds with probability at least $1-\delta$, if $\bTheta$ is $(1/2, \delta/n, 1)$ oblivious subspace embedding. 

Furthermore, let us define rounding error vectors 
$$\bDelta \bq'_i := \bhq'_i -(\bhw_i-\bhQ_{i-1}[\bhR]_{(1:i-1,i)}) \text{ and } \bDelta \bq_i := \bhq_i - \bhq_i'/\hat{r}_{i,i}$$ in Steps 3 and 6 of~\cref{alg:RGM}. Then, the standard worse-case scenario rounding analysis yields~\cref{eq:numepsembedding3stand}  {(for derivation, see~\cref{skroundofferrors})}. Following the arguments from~\cref{skroundofferrors}, we assume that $\bTheta$ satisfies~\cref{eq:numepsembedding3}. It follows from~\cref{thm:thetadeltax2} and the union bound argument that this property holds under the {probabilistic} rounding model with probability at least $1-4\delta$, if $\bTheta$ is  $(1/8, m^{-1}\binom{n}{d}^{-1}\delta, d)$, {with $d = \mathcal{O}(\log(m/\delta))$,} oblivious subspace embedding. By using the bounds~\cref{eq:ell_2bounds} we conclude that~\cref{thm:Thetaasmpts} hold under the {probabilistic} rounding model with probability at least $1-\delta$, if $\bTheta$ is a Rademacher matrix with {$k = \mathcal{O}(\log(n)\log(m/\delta))$} rows or SRHT with $k = \mathcal{O}(\log^2(n)\log^2(m/\delta))$ rows. Note that these properties should also hold under many other (possibly deterministic) rounding models.  The classical worse-case scenario model, however, entails $D = \sqrt{1+\varepsilon} \sqrt{n}$ in~\cref{eq:numepsembedding3}. The numerical stability bounds in this case can be deduced from~\cref{thm:maintheorem1,thm:maintheorem2} by letting $u_{crs} =  u_{crs} \sqrt{n}$ and $u_{fine} = u_{fine} \sqrt{n}$.

\begin{assumptions} \label{thm:Thetaasmpts}
	It is assumed that
	\begin{equation} \label{eq:Thetanorm}
	\| \bTheta \|_\Frob \leq \sqrt{1+\varepsilon} \sqrt{n}, 
	\end{equation}
	with $\varepsilon \leq 1/2$. Furthermore, in~\cref{alg:RGM} we assume that
	\begin{subequations} \label{eq:numepsembedding3stand} 
		\begin{align}
		| \bDelta \bq'_i| &  \leq 1.02 u_{crs} (| \bhw_i | + i | \bhQ_{i-1}| |[\bhR]_{(1:i-1,i)} |), \\
		| \bDelta \bq_i| &\leq   { u_{fine}} |\bhq_i'/\hat{r}_{i,i}|,
		\end{align}
	\end{subequations}
	and
	\begin{subequations} \label{eq:numepsembedding3} 
		\begin{align}
		\|\bTheta \bDelta \bq'_i\| &\leq 1.02 u_{crs} D \| |\bhw_i| + i | \bhQ_{i-1}| |[\bhR]_{(1:i-1,i)}| \| , \\
		\| \bTheta \bDelta \bq_i\| &\leq {u_{fine}} D \|\bhq_i'/\hat{r}_{i,i} \|,
		\end{align}
	\end{subequations}
	with $D = \sqrt{1+\varepsilon}$, $\varepsilon  \leq 1/2$, $1 \leq i \leq m$.
\end{assumptions}

\subsection{Stability analysis} \label{stability}

The results in this subsection shall rely on the condition that $\bTheta$ satisfies the $\varepsilon$-embedding property 
for $\bhQ$ and $\bhW$.
A priori analysis to satisfy this property with (high) user-specified probability of success is provided in~\cref{l2embedding}. Furthermore, in~\cref{l2embedding} we also provide a way to efficiently certify that  $\bTheta$ is {an} $\varepsilon$-embedding for $\bhQ$ and $\bhW$.	
Then the stability of~\cref{alg:RGM} can be characterized by coefficients
$$ \Delta_m = \|\bI_{m \times m} - \bhS_m^\mathrm{T} \bhS_m \|_\Frob  \text{ and } \tilde{\Delta}_m = \|\bhP_m - \bhS_m \bhR_m\|_\Frob/\|\bhP_m\|_\Frob,$$
as it is shown in~\cref{thm:maintheorem1}.

\begin{theorem}\label{thm:maintheorem1}
	{Assume that} 
	$$ {100} m^{1/2} n^{3/2} u_{fine} \leq u_{crs} \leq 0.01 m^{-1}, \text{ and } n\geq 100,$$
	{along with~\cref{thm:Thetaasmpts}.}
	{If $\bTheta$ is {an} $\varepsilon$-embedding for $\bhQ$ and $\bhW$ from \Cref{alg:RGM}, with $\varepsilon \leq 1/2$, and if $\Delta_m,\tilde{\Delta}_m \leq 0.1$,}
	then the following inequalities hold:
	$$ (1+\varepsilon)^{-1/2}(1-{\Delta_m}- 0.1 u_{crs}) \leq \sigma_{min}(\bhQ)  \leq \sigma_{max}(\bhQ) \leq (1-\varepsilon)^{-1/2}(1+{\Delta_m}+ 0.1 u_{crs})$$
	and 
	$$\|\bhW - \bhQ \bhR\|_\Frob \leq  3.7 u_{crs} m^{3/2} \|\bhW\|_\Frob. $$	
	\begin{proof}
		See supplementary material.
	\end{proof}
\end{theorem}

According to~\cref{thm:maintheorem1}, the numerical stability of~\cref{alg:RGM} can be ensured by guaranteeing that $\Delta_m$ and $\tilde{\Delta}_m$ are sufficiently small {(along with the $\varepsilon$-embedding property of $\bTheta$)}. This can be done by employing a sufficiently accurate backward-stable solver to the least-squares problem in Step 2.
Below we provide theoretical bounds for $\Delta_m$ and $\tilde{\Delta}_m$ in this case.
\begin{theorem}\label{thm:maintheorem2}
	{Consider~\cref{alg:RGM} utilizing QR factorization based on Householder transformation or Givens rotations for computing the solution to the  least-squares problem in Step 2. }
	
	{Under~\cref{thm:Thetaasmpts},} if $\bTheta$ is {an} $\varepsilon$-embedding for $\bhQ_{m-1}$ and $\bhW$, with $\varepsilon \leq 1/2$, and if
	\begin{align*} 
	u_{crs} &\leq 10^{-3} \cond(\bhW)^{-1} m^{-2}, \\
	u_{fine}   &\leq (100 m^{1/2} n^{3/2} + 10^4 m^{3/2} k)^{-1}u_{crs},
	\end{align*}
	then $\Delta_m$ and $\tilde{\Delta}_m$ are bounded by
	\begin{align} 
	\tilde{\Delta}_m 
	& \leq 4.2 u_{crs}  m^{3/2} \|\bhW\|_\Frob/\|\bhP\|_\Frob \leq 6 u_{crs} m^{3/2}, \label{eq:maintheorem21} \\
	\Delta_m 
	&\leq 20 u_{crs}   m^{2} \cond{(\bhW)}. \label{eq:maintheorem22}
	\end{align}
	\vspace*{-0.5cm}
	\begin{proof}
		See supplementary material.
	\end{proof}
\end{theorem}

\begin{remark}
	In general, the result of~\cref{thm:maintheorem2} holds for any least-squares solver {in Step 2} as long as the following backward-stability property is satisfied: 	
	\begin{gather*}
	[\bhR]_{(1:i-1,i)}  = \arg \min_{\by} \| (\bhS_{i-1}+ \bDelta\bS_{i-1})\by - (\bhp_i + \bDelta\bp_i) \|,\text{with} \\
	\| \bDelta\bS_{i-1} \|_\Frob \leq 0.01 u_{crs}   \| \bhS_{i-1} \|,~~   \| \bDelta\bp_i \| \leq 0.01 u_{crs}   \|\bhp_i \|. 
	\end{gather*}
\end{remark}

\begin{remark} \label{rmk:proofmaintheorem2}
	{Notice that \Cref{thm:maintheorem2} requires $\bTheta$ to be {an} $\varepsilon$-embedding for $\bhQ_{m-1}$ and not $\bhQ$. This observation will become handy  for proving the $\varepsilon$-embedding property for $\bhQ$ in~\cref{l2embedding} by using induction on $m$.}
\end{remark}
\Cref{thm:maintheorem1,thm:maintheorem2} imply a stable QR factorization for {working} unit roundoff $u_{crs}$ independent of {the high dimension $n$.}

In some cases, obtaining a priori guarantees with \Cref{thm:maintheorem2} can be an impractical task due to the need to estimate $\cond(\bhW)$. Furthermore, one may want to use~\cref{alg:RGM} with a higher value of $u_{crs}$ than is assumed in~\cref{thm:maintheorem2}, possibly with a bigger gap between the values of $u_{crs}$ and $u_{fine}$. In such cases, the computed QR factorization can be efficiently certified a posteriori by computing $\Delta_m$ and $\tilde{\Delta}_m$ and using~\cref{thm:maintheorem1} with no operations on high-dimensional vectors and matrices, and the estimation of $\mathrm{cond}(\bhW)$.

\subsection{Epsilon embedding property} \label{l2embedding}
The stability analysis in~\cref{stability} holds if $\bTheta$ satisfies the $\varepsilon$-embedding property for $\bhQ$ and $\bhW$. In this section we provide a priori and a posteriori analysis of this property.

\subsubsection*{A priori analysis} \label{apriori}

Let us consider the case when $\bhW$ and $\bTheta$ are independent of each other. Then it follows directly from~\cref{def:oblemb} that, if $\bTheta$ is $(\varepsilon, \delta, m)$ oblivious $\ell_2$-subspace embedding, then it satisfies the $\varepsilon$-embedding property for $\bhW$ with high probability. Below, we provide a guarantee that{,} in this case $\bTheta$ will also satisfy {an} $\varepsilon$-embedding property for $\bhQ$ with moderately increased value of $\varepsilon$.
\begin{proposition}\label{thm:aprioribound}
	Consider~\cref{alg:RGM} {using the Givens or Householder least-squares solver in Step 2}, and computed with unit roundoffs 
	\begin{align*} 
	u_{crs} &\leq 10^{-3} \cond(\bhW)^{-1} m^{-2}, \\
	u_{fine}   &\leq (100 m^{1/2} n^{3/2} + 10^4 m^{3/2} k)^{-1}u_{crs}.
	\end{align*}
	{Under~\cref{thm:Thetaasmpts},} if $\bTheta$ is {an} $\varepsilon$-embedding for $\bhW$, with $\varepsilon \leq 1/4$, then it satisfies the $\varepsilon'$-embedding property for $\bhQ$, with 
	$$\varepsilon' = 2 \varepsilon + 180 u_{crs} m^2 \cond(\bhW). $$
	\vspace*{-0.5cm}	
	\begin{proof}  
		See supplementary material.
	\end{proof}
\end{proposition}

When $\bhW$ is generated depending on $\bTheta$ (as we have in Arnoldi process in~\cref{application}), the a priori analysis for the $\varepsilon$-embedding property for  $\bhW$ can be nontrivial and pessimistic. Nevertheless,  $\bTheta$ can still be expected to be {an} $\varepsilon$-embedding {because there is only} a minor correlation of the rounding errors with $\bTheta$. When there is no a priori  guarantee on the quality of $\bTheta$ or the guarantee is pessimistic, it can be important to be able to certify the $\varepsilon$-embedding property a posteriori, which is discussed next.

\subsubsection*{A posteriori certification} 
The quality of $\bTheta$ can be certified by providing an upper bound $\bar{\omega}$ for the minimum value $\omega$ of $\varepsilon$, for which $\bTheta$ satisfies the $\varepsilon$-embedding property for $\bV$. The matrix $\bV$ can be chosen as $\bhQ$ or $\bhW$. 
The considered a posteriori bound $\bar{\omega}$ is probabilistic. 
We proceed by introducing an (additional to $\bTheta$) sketching matrix $\bPhi$ used solely for the certification so that it is randomly independent of $\bV$ and $\bTheta$. {For efficiency, this matrix should be of size no more than the size of $\bTheta$. In practice, an easy and robust way is to use $\bPhi$ and $\bTheta$ of same size.} Define parameters $\varepsilon^*$ and $\delta^*$ characterizing respectively, the accuracy of $\bar{\omega}$ (i.e., its closeness to $\omega$) and the probability of failure for $\bar{\omega}$ to be an upper bound. Then we can use the following results from~\parencite{balabanov2019randomized2}(see~\cref{thm:cert01,prop:cert2}).
\begin{proposition}[Corollary of~Proposition 5.3 in \parencite{balabanov2019randomized2}] \label{thm:cert01}
	Assume that $\bPhi$ is a $(\varepsilon^*, \delta^*, 1)$-oblivious subspace embedding.
	Let $\bV = \bhQ$ or $\bhW$, $\bV^\bTheta = \bTheta \bV$, and $\bV^\bPhi = \bPhi \bV$. 	
	Let $\bX$ be a matrix such that $\bV^\bPhi\bX$ is orthonormal. If, 
	$$\bar{\omega} =  \max \{1-(1-\varepsilon^*) \sigma^2_{min}(\bV^\bTheta\bX),  (1+\varepsilon^*) \sigma^2_{max}(\bV^\bTheta\bX)-1  \} <1,   $$
	then $\bTheta$ is a $\bar{\omega}$-embedding for $\bV$,  with probability at least $1-\delta^*$.	
\end{proposition}
\begin{proposition}[Corollary of~Proposition 5.4 in \parencite{balabanov2019randomized2}]\label{prop:cert2}
	In~\cref{thm:cert01}, if $\bPhi$ is a $\varepsilon'$-embedding for $\bV$, then $\bar{\omega}$ satisfies $$\bar{\omega} \leq (1+\varepsilon^*)(1-\varepsilon')^{-1}(1+\omega) -1.$$	
\end{proposition}
First we notice that the coefficient $\bar{\omega} $ in \Cref{thm:cert01} can be efficiently computed from the two sketches of $\bV$ with no operations on high-dimensional vectors and matrices. For efficiency, (at iteration $i$) the sketches $\bPhi \bhq'_i$ and $\bPhi \bhw_i$ may be computed along with, respectively, $\bhs'_i$ in Step 4 and $\bhp_i$ in Step 1 of~\cref{alg:RGM}. The matrix $\bX$ can be obtained (possibly in implicit form, e.g., as inverse of upper-triangular matrix) {with standard orthogonal decomposition algorithms such as the QR factorization or the singular value decomposition} performed in sufficient precision.   It follows that $\bar{\omega}$ is an upper bound for $\omega$ with probability at least $1-\delta^*$, if $\bPhi$ is a $(\varepsilon^*, \delta^*, 1)$-oblivious subspace embedding. This property of $\bPhi$ has to be guaranteed a priori, e.g, from the theoretical bounds~\cref{eq:ell_2bounds} or~\parencite[Lemma 5.1]{achlioptas2003database}. For instance, it is guaranteed to hold for Rademacher matrices with $k \geq 2({\varepsilon^*}^2/2-{\varepsilon^*}^3/3)^{-1} \log(\delta^*/2)$ rows, which in particular becomes $k > 530$ for $\varepsilon^* =1/4$ and $\delta^* = 0.1\%$ and any $m$ and $n$.

In~\cref{prop:cert2}, the closeness of $\bar{\omega}$ and $\omega$ is guaranteed if $\bPhi$ is a $\varepsilon'$-embedding for $\bV$ (for some given $\varepsilon'$). This condition shall be satisfied with probability at least $1-\delta'$ (for some given $\delta'$), if $\bPhi$ is an $(\varepsilon', \delta', m)$ oblivious $\ell_2$-subspace embedding. This fact is not required to be guaranteed a priori, which allows to choose the size for $\bPhi$ (and $\bTheta$) based on practical experience and still have a certification.  

In practice, the random projections $\bV^\bTheta = \bTheta \bV$ and $\bV^\bPhi = \bPhi \bV$ can be computed only approximately due to rounding errors. In such case, it can be important to provide a stability guarantee for the computed value of $\bar{\omega}$ ({given in}~\cref{prop:cert1}). {In~\cref{prop:cert1}, along with \cref{thm:Thetaasmpts} for $\bTheta$, we also assume} that
\begin{equation} \label{eq:Phinorm}
\| \bPhi \|_\Frob \leq \sqrt{1+\varepsilon} \sqrt{n},~\text{ and } \|\bV^\bPhi \|_\Frob \geq \sqrt{1-\varepsilon} \| \bV \|_\Frob,
\end{equation}
with $\varepsilon \leq 1/2$. These properties hold with probability at least $1-2\delta^*$, if $\bPhi$ is $(1/2, \delta^*/n, 1)$ oblivious subspace embedding. 
\begin{proposition} \label{prop:cert1}
	Let $\bV = \bhQ$ or $\bhW$, $\bhV^\bTheta = \fl(\bTheta \cdot \bV)$, and $\bhV^\bPhi = \fl(\bPhi \cdot \bV)$. 
	Assume that the  sketches are computed with unit roundoff $u_{fine}$ satisfying
	$$  100 n^{3/2} m^{1/2} u_{fine} \leq u_{crs} \leq {\cond(\bhV^\bPhi)^{-1}}.$$ 
	{Assume that $\bPhi$ satisfies~\cref{eq:Phinorm} and that $\bTheta$ satisfies~\cref{thm:Thetaasmpts}.}
	Let $\bhX$ be a matrix such that $\bhV^\bPhi\bhX$ is orthonormal. Define
	$$\hat{\bar{\omega}} =  \max \{1-(1-\varepsilon^*) \sigma^2_{min}(\bhV^\bTheta\bhX),  (1+\varepsilon^*) \sigma^2_{max}(\bhV^\bTheta\bhX)-1  \}.$$
	Then we have if $\bar{\omega} \leq 1$:
	$$ | \bar{\omega} -  \hat{\bar{\omega}}| \leq u_{crs}\cond(\bhV^\bPhi).$$
	\begin{proof}
		See supplementary material.
	\end{proof}
\end{proposition}
It follows form~\cref{prop:cert1} that the computed value of $\bar{\omega}$ is approximately equal to the exact one if $u_{crs}\cond(\bhV^\bPhi)$ is sufficiently small. This implies the following computable certificate for the quality of $\bTheta$:
$$ \omega \leq \bar{\omega} \leq \hat{\bar{\omega}} + u_{crs}\cond(\bhV^\bPhi), $$
which holds with probability at least $1-\mathcal{O}(\delta^*)$.
Furthermore, \Cref{prop:cert1} is consistent: $\cond(\bhV^\bPhi)$ is guaranteed to be sufficiently small, namely $\mathcal{O}(\cond(\bV))$, if $\bPhi$ is a $\varepsilon$-embedding for $\bV$.

\begin{remark}[Randomized Gram-Schmidt algorithm with multiple sketches] \label{rmk:multRGS}
	The certification of $\bTheta$ can be performed  at each iteration of~\cref{alg:RGM} (by letting $m = i$ in the above procedure). In this way, one can be able to detect the iteration $i$ (if there is any) with not sufficient quality of $\bTheta$ and switch to new sketching matrix, which is randomly independent from $\bV$. This allows to make sure that the used $\bTheta$ is {an} $\varepsilon$-embedding for $\bV$ at each iteration. 
	We leave the development of randomized GS algorithm with multiple sketches for future research.  
\end{remark}

\section{Application to Arnoldi process and GMRES} \label{application}
In this section we employ the randomized GS process to solving  high-dimensional non-singular (possibly non-symmetric) systems of equations of the form
\begin{equation} \label{eq:insystem}
\bA \bx = \bb,
\end{equation}
with GMRES method. 
Without loss of generality we here assume that~\cref{eq:insystem} is normalized so that $\|\bb \| = \|\bA\| = 1$.   

An order-$j$ Krylov subspace is defined as 
$$\mathcal{K}_{j}(\bA, \bb) := \mathrm{span} \{\bb, \bA \bb, \hdots, \bA^{j-1} \bb \}.$$
The GMRES method consists in approximation of $\bx$ with a projection $\bx_{m-1}$ in $\mathcal{K}_{m-1}(\bA, \bb)$ that
minimizes the residual norm 
$$\| \bA \bx_{m-1} - \bb \|. $$
To obtain a projection $\bx_{m-1}$, the GMRES method first proceeds with constructing orthonormal basis of $\mathcal{K}_m(\bA, \bb)$ with Arnoldi process (usually based on GS orthogonalization). Then the coordinates of $\bx_{m-1}$ in the Arnoldi basis are found by solving a (small) transformed least-squares problem.

\subsection{Randomized GS-Arnoldi process} 
The Arnoldi basis can be constructed recursively by taking the first basis vector $\bq_1$ as normalized right-hand-side vector $\bb$ and each new vector $\bq_{i+1}$ as $\bA \bq_{i}$ orthonormalized against the previously computed basis $\bq_{1}, \hdots, \bq_{i}$. This procedure then produces orthonormal matrix $\bQ_m$ satisfying the Arnoldi identity 
$$ \bA \bQ_{m-1} = \bQ_{m} \bH_m, $$
where $\bH_m$ is upper {Hessenberg} matrix. The Arnoldi algorithm can be viewed as a column-oriented QR factorization of matrix $[\bb, \bA \bQ_{m-1}]$. In this case, the R factor $\bR_m$ and the {Hessenberg} matrix $\bH_m$ satisfy the relation $\bH_m = [\bR_m]_{(1:m, 2:m)}.$

Below, we propose a randomized Arnoldi process based on randomized GS algorithm from~\cref{RGS} for computing {the} Krylov basis orthonormal with respect to the sketched product $\langle \bTheta \cdot, \bTheta \cdot \rangle$, rather than $\ell_2$-inner product as in standard methods (see~\cref{alg:RGS-Arnoldi}). {This process will serve as the core for the randomized GMRES method in~\cref{RGMRES} }.

\begin{algorithm}[h] \caption{RGS-Arnoldi algorithm} \label{alg:RGS-Arnoldi}
	\begin{algorithmic}
		\STATE{\textbf{Given:} $n \times n$ matrix $\bA$, $n \times 1$ vector $\bb$, $k \times n$ matrix $\bTheta$ with $k \ll n$, parameter $m$.}
		\STATE{\textbf{Output}: $n \times m$ factor $\bQ_m$ and $m \times m$ upper triangular factor $\bR_m$.}
		\STATE{1. Set $\bw_1 = \bb$. }
		\STATE{2. Perform $1$-st iteration of~\cref{alg:RGM}.}
		\FOR{ $i=2:m$} 		
		\STATE{3. Compute $\bw_i = \bA \bq_{i-1}$. }\COMMENT{macheps: $u_{fine}$}
		\STATE{4. Perform $i$-th iteration of~\cref{alg:RGM}.}
		\ENDFOR
		\STATE{5. (Optional) compute $\Delta_m$ and $\tilde{\Delta}_m$. \\~~~ Use~\cref{thm:maincorollary1} to certify the output.}\COMMENT{macheps: $u_{fine}$}
	\end{algorithmic}
\end{algorithm}

In~\cref{alg:RGS-Arnoldi}, the computation of the matrix-vector product in Step~3 with the fine unit roundoff $u_{fine}$ is assumed to have only a minor impact on the overall computational costs. This can be the case, for instance, when {the} matrix $\bA$ is sparse or structured. {Furthermore, if needed, the matrix-vector product can be computed also with a larger unit roundoff as long as the associated error satisfies
	$\|\bhw_i - \bA \bhq_{i-1}\| = \mathcal{O}(u_{crs} m^{-1/2}) \|\bhq_{i-1}\| $ required by~\cref{thm:maincorollary1,thm:maincorollary2}.}

Let us now address the accuracy of~\cref{alg:RGS-Arnoldi}. Clearly, if $\bTheta$ is {an} $\varepsilon$-embedding for $\mathcal{K}_m(\bA, \bb)$, then~\cref{alg:RGS-Arnoldi} in infinite precision arithmetic produces  a well-conditioned basis matrix $\bQ_{m}$ satisfying the Arnoldi identity. In addition, we clearly have $\mathrm{range}{(\bQ_{m})} = \mathrm{range}{(\bW_{m})} =\mathcal{K}_m(\bA, \bb)$. Since $\mathcal{K}_m(\bA, \bb)$ and $\bTheta$ are independent, the matrix $\bTheta$ can be readily chosen as $(\varepsilon,\delta,m)$ oblivious $\ell_2$-subspace embedding to have the $\varepsilon$-embedding property  with high probability.

Numerical stability of~\cref{alg:RGS-Arnoldi} in finite precision arithmetic can be derived directly from~\cref{thm:maintheorem1,thm:maintheorem2} characterizing the stability of the randomized GS algorithm.  
\begin{proposition}\label{thm:maincorollary1}
	Assume that 
	$$ {100 }   m^{1/2} n^{3/2} u_{fine}\leq u_{crs} \leq 0.01 m^{-1}, \text{ and } n\geq 100,$$
	{along with \cref{thm:Thetaasmpts}.} If $\bTheta$ satisfies the $\varepsilon$-embedding property for $\bhQ_m$ and $\bhW_m$ with $\varepsilon \leq 1/2$ and $\Delta_m,\tilde{\Delta}_m \leq 0.1$, then we have
	$$  (1+\varepsilon)^{-1/2} (1- \Delta_m - 0.1 u_{crs})  \leq \sigma_{min}(\bhQ_m) \leq \sigma_{max}(\bhQ_m) \leq  (1-\varepsilon)^{-1/2} (1+ \Delta_m + 0.1 u_{crs}). $$
	We also have,
	$$ (\bA+\bDelta\bA) \bhQ_{m-1} = \bhQ_{m} \bhH_{m}$$
	for some matrix $\bDelta\bA$ with $\mathrm{rank}(\bDelta\bA)<m$  and $\|\bDelta\bA\|_\Frob \leq 15 u_{crs} m^2.$ 
	\begin{proof}
		See supplementary material.
	\end{proof}
\end{proposition}
It follows from \Cref{thm:maincorollary1} that the stability of the proposed RGS-Arnoldi algorithm can be guaranteed by computing or bounding a priori coefficients $\Delta_m$ and $\tilde{\Delta}_m$. \Cref{thm:maincorollary1} can be viewed as a backward stability characterization, since it implies that
$$ \mathrm{range}(\bhQ_{m-1}) = \mathcal{K}_{m-1}(\bA + \bDelta \bA , \bb +\bDelta \bb),$$
where $\|\bDelta\bA\|_\Frob \leq 15 u_{crs} m^{2}$ and $\| \bDelta \bb \| = \|\hat{r}_{1,1} \bhq_1 - \bb\| \leq u_{fine} $. In other words, according to~\cref{thm:maincorollary1}, the output of~\cref{alg:RGS-Arnoldi} in finite precision arithmetic is guaranteed to produce a well-conditioned basis $\bhQ_{m-1}$ for the Krylov space of a slightly perturbed matrix $\bA$ and vector $\bb$.

Let us next provide a priori bounds for $\Delta_m$ and $ \tilde{\Delta}_m$.  Define parameter $$\tau(\bhQ_{m-1}) = \min_{\by_{m-1} \in \mathrm{range}(\bhQ_{m-1})} \|  \bA \by_{m-1} - \bb \|, $$
representing the best attainable residual error with the computed Krylov basis.

\begin{proposition}\label{thm:maincorollary2}
	Consider~\cref{alg:RGS-Arnoldi} {using the Givens or Householder least-squares solver in Step 2 of~\cref{alg:RGM},  and}
	\begin{align*} 
	u_{crs} &\leq 10^{-4} \tau(\bhQ_{m-1}) \cond({\bA})^{-1} m^{-2}, \\
	u_{fine}   &\leq (100 m^{1/2} n^{3/2} + 10^4 m^{3/2} k)^{-1}u_{crs}.
	\end{align*}
	{Under~\cref{thm:Thetaasmpts},} if $\bTheta$ satisfies the $\varepsilon$-embedding property for $\bhQ_{m-1}$ and $\bhW_m$ with $\varepsilon \leq 1/2$, then $\Delta_m$ and $\tilde{\Delta}_m$ are bounded by
	\begin{align*} 
	\tilde{\Delta}_m & \leq 6 u_{crs} m^{3/2}, \\
	\Delta_m &\leq 160 u_{crs}   m^{2} \cond{(\bA)} \tau(\bhQ_{m-1})^{-1}. 
	\end{align*}
	\vspace*{-0.5cm}	
	\begin{proof}
		See supplementary material.
	\end{proof}
\end{proposition}
\Cref{thm:maincorollary2} guarantees numerical stability of~\cref{alg:RGS-Arnoldi}, if $$\tau(\bhQ_{m-1})  \geq u_{crs} P(m) \cond{(\bA)},$$ where $P(m) = \mathcal{O}(m^2)$ is some low-degree polynomial.  Clearly, if $\tau_0$ is the desired tolerance for the GMRES solution, then one is required to use unit roundoff $u_{crs} \leq \tau_0 P(m)^{-1} \cond{(\bA)}^{-1}$. 

Both~\cref{thm:maincorollary1,thm:maincorollary2} hold if $\bTheta$ satisfies the $\varepsilon$-embedding property for $\bhQ_m$ and $\bhW_m$ with $\varepsilon \leq 1/2$. This assumption comes naturally from probabilistic characteristics of rounding errors and oblivious embeddings. In particular, we can think of similar considerations as in~\cref{skroundofferrors} to justify the $\varepsilon$-embedding property of $\bTheta$ for  $\mathcal{K}_{m}(\bA + \bDelta \bA , \bb +\bDelta \bb)$, where matrix $\bDelta \bA$ has rows with entries that are independent centered random variables. The a priori analysis of the $\varepsilon$-embedding property for a perturbed Krylov space, however, is not as trivial, and is left for future research. Note that the output of~\cref{alg:RGS-Arnoldi} can be proven reliable a posteriori by efficient certification of the $\varepsilon$-embedding property  with the procedure from~\cref{l2embedding}.   

\subsection{Randomized GMRES} \label{RGMRES}

Randomized GMRES method is directly derived from the randomized Arnoldi iteration. Let $\bhQ_m$ and $\bhH_m$, be the basis matrix and the Hessenberg matrix computed with~\cref{alg:RGS-Arnoldi}. (Randomized) GMRES method then proceeds with obtaining {the} solution $\by_{m-1}$ to the following small least-squares problem  
\begin{equation} \label{eq:hls}
\by_{m-1} = \arg \min_{\bz_{m-1} \in \mathbb{R}^{m-1}} \| \bhH_m \bz_{m-1} - \hat{r}_{1,1} \be_1 \|,
\end{equation}
yielding an approximate solution $\bx_{m-1} = \bhQ_{m-1} \by_{m-1}$ to~\cref{eq:insystem}. The underlined least-squares problem can be efficiently solved with a QR factorization based on Givens rotations or Householder transformation (or any other methods) in sufficient precision.  

{In infinite precision arithmetic, the orthogonality of $\bQ_m$ implies that  solving~\cref{eq:hls} is equivalent to minimizing \emph{sketched} norm of the residual. More specifically, we have
	\begin{equation} \label{eq:hls2}
	\begin{split}
	\min_{\bz_{m-1} \in \mathbb{R}^{m-1}} \| \bH_m \bz_{m-1} - r_{1,1} \be_1 \| &= \min_{\bz_{m-1} \in \mathbb{R}^{m-1}} \| \bTheta \bQ_m (\bH_m \bz_{m-1} - r_{1,1} \be_1) \| \\
	&= \min_{\bz_{m-1} \in \mathbb{R}^{m-1}} \| \bTheta (\bA \bQ_{m-1} \bz_{m-1} - \bb) \|.
	\end{split}
	\end{equation}
	Thus, the GMRES solution $\bx_{m-1}$ minimizes the residual error up to a factor  $\sqrt{\frac{1+\varepsilon}{1-\varepsilon}}$, provided $\bTheta$ is an $\varepsilon$-embedding for $\bQ_m$. }

Numerical stability of the randomized GMRES can be characterized by using~\cref{thm:maincorollary1} {that yields the following result.}
\begin{proposition} \label{thm:GMRESstab}
	{Let~\cref{thm:Thetaasmpts} hold. }Assume that $\bTheta$ is an $\varepsilon$-embedding for $\bhQ_m$ and $\bhW_m$, with $\varepsilon \leq 1/2$ and $\Delta_m,\tilde{\Delta}_m \leq 0.1$, then we have
	$$ \| (\bA+\bDelta \bA) \bx_{m-1} - (\bb +\bDelta \bb) \| \leq \cond(\bhQ_m) \min_{\bv \in \mathcal{K}_{m-1}(\bA + \bDelta \bA , \bb +\bDelta \bb)} \| (\bA+\bDelta \bA) \bv - (\bb +\bDelta \bb) \|, $$
	for some matrix $\bDelta\bA$ and vector $\bDelta \bb$ with $\|\bDelta\bA\|_\Frob \leq 15 u_{crs} m^{2}$ and $\| \bDelta \bb \| \leq u_{fine}$. 
	\begin{proof}
		See supplementary material.
	\end{proof}
\end{proposition}
Notice that for sufficiently small $\Delta_m$ and $\tilde{\Delta}_m$, $\cond(\bhQ_m)$ is close to $\sqrt{\frac{1+\varepsilon}{1-\varepsilon}}$. Consequently,~\cref{thm:maincorollary1} guarantees that $\bx_{m-1}$ is a quasi-optimal minimizer of the residual {error} over a (slightly) perturbed Krylov space.  

\section{Numerical experiments} \label{numexp}

In this section the proposed methodology is verified in a series of numerical experiments and compared against classical methods. In the randomized algorithms, several sizes $k$ and distributions for the sketching matrices $\bTheta$ and $\bPhi$ are tested. We use in Step 2 of~\cref{alg:RGM} the Householder least-squares solver.   There was not detected any significant difference in performance (i.e., stability or accuracy of approximation for the same $k$) between Rademacher (or Gaussian) distribution  and (P-)SRHT, even though the theoretical bounds for (P-)SRHT are worse. Therefore, in this section we present only the results for the (P-)SRHT distribution. 

For better presentation, the orthogonality of the sketch $\bS_m$ is here measured by the condition number $\cond(\bS_m)$ instead of the coefficient $\Delta_m  = \|\bI_{m \times m} - \bS_m^\mathrm{T} \bS_m \|_\Frob$ as in the previous sections. 


\subsection{Construction of an orthogonal basis for synthetic functions} \label{Ex1}
Let us first consider construction of an orthogonal basis approximating the functions:
\begin{equation*}
f_\mu(x) = \frac{\sin\left (10(\mu+x) \right)}{\cos \left(100(\mu-x) \right)+1.1},~x \in [0,1],
\end{equation*}
for parameter values $\mu \in [0,1]$. 

The function's domain is discretized with $n = 10^6$ evenly spaced points $x_j$, while the parameter set is discretized with $m = 300$ evenly spaced points $\mu_j$. Then a QR factorization of the matrix $[\bW]_{i,j} = f_{\mu_j}(x_i)$, $1\leq i \leq n$, $1 \leq j \leq m$, is performed with standard versions (CGS, MGS, and CGS2) of Gram-Schmidt process, along with the randomized version of the process, given by~\cref{alg:RGM}. 
The classical algorithms are executed in float32 format with unit roundoff $\approx 10^{-8}$. \Cref{alg:RGM} is first executed using a unique float32 format for all the arithmetic operations, i.e., by taking $u_{crs} = u_{fine} \approx 10^{-8}$. Then, the results are compared to  \Cref{alg:RGM} under the multi-precision model executing Step 3 in float32, while executing other operations in float64, i.e., by taking $u_{crs} \approx 10^{-8}$ and $u_{fine} = 10^{-16}$. Note that the execution of~\cref{alg:RGM} with the unique float32 format has nearly the same computational cost as with the mixed float32/float64 formats. Furthermore, as was argued in~\cref{peranal}, the RGS algorithm here requires twice less flops\footnote{We did not take into consideration the flops associated with the solutions of $k \times (i-1)$ least-squares problems in Step 2 of~\cref{alg:RGM}, which will become irrelevant for larger systems. They could be reduced by solving the least-squares problems (iteratively) with normal equation.} and {data passes} than CGS and, respectively, four times less flops and {data passes} than CGS2. Moreover,  unlike MGS, it is implemented by using BLAS-2 routines for standard high-dimensional operations. 

\begin{figure}[!h]
	\centering
	\begin{subfigure}{.35\textwidth}
		\centering  
		\includegraphics[width=\textwidth]{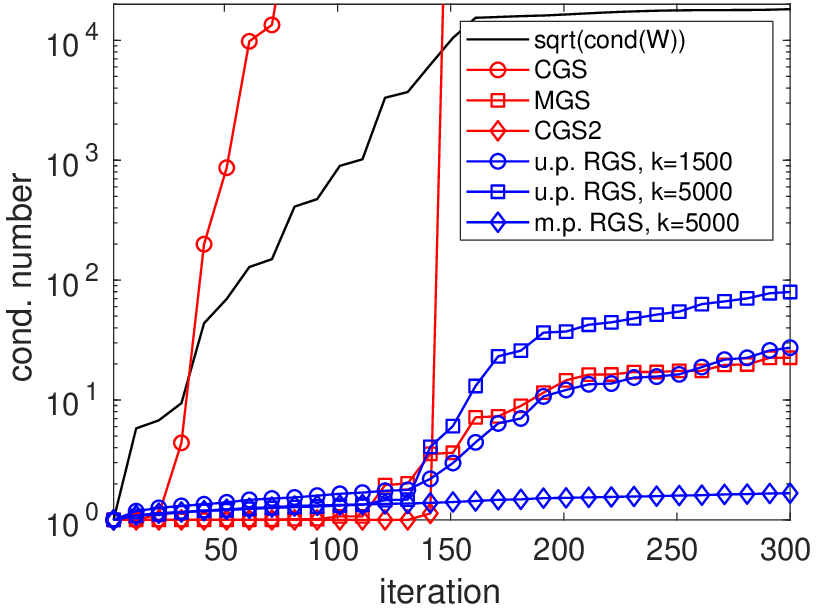}
		\caption{Cond. number of $\bQ_i$}
		\label{fig:Ex1_1a}
	\end{subfigure} \hspace{.1\textwidth}
	\begin{subfigure}{.35\textwidth}
		\centering
		\includegraphics[width=\textwidth]{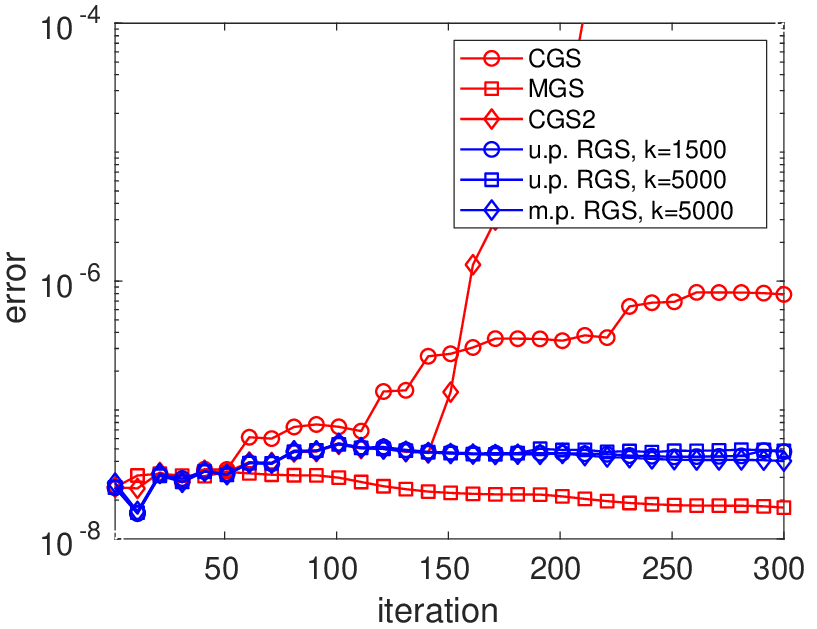}
		\caption{Error $\|\bW_i - \bQ_i \bR_i\|/\|\bW_i\|$}
		\label{fig:Ex1_1b}
	\end{subfigure} \\
	\begin{subfigure}{.35\textwidth}
		\centering  
		\includegraphics[width=\textwidth]{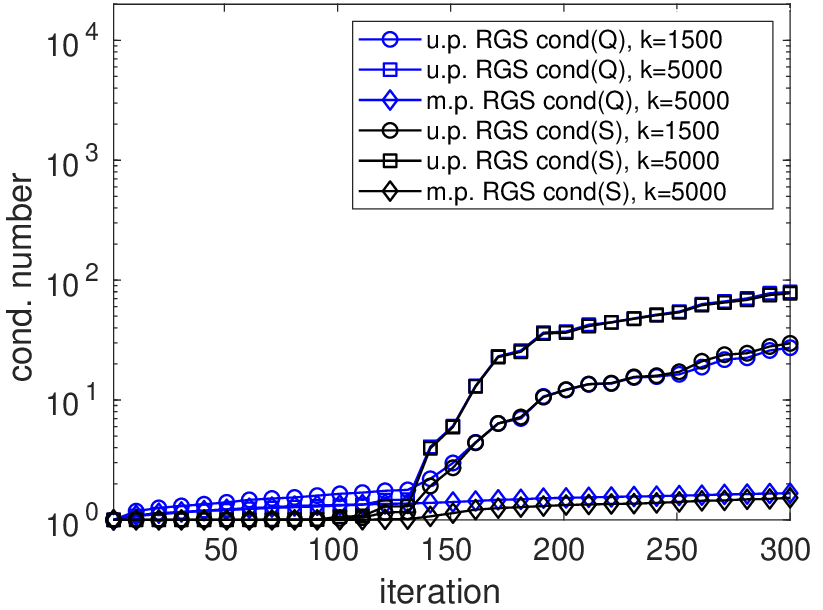}
		\caption{Cond. numbers of $\bQ_i$ and $\bS_i$.}
		\label{fig:Ex1_1c}
	\end{subfigure} \hspace{.1\textwidth}
	\begin{subfigure}{.35\textwidth}
		\centering
		\includegraphics[width=\textwidth]{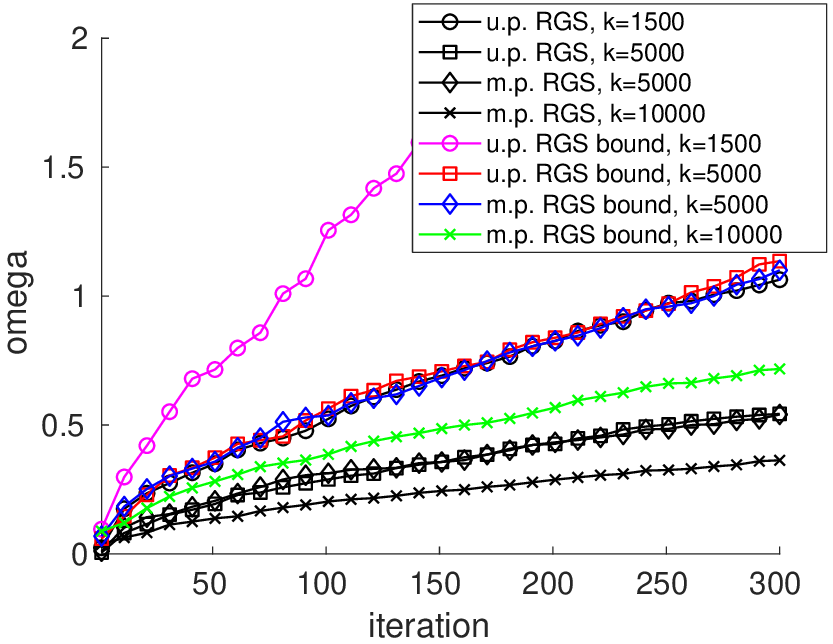}
		\caption{Coeff. $\omega$ and its upper bound}
		\label{fig:Ex1_1d}
	\end{subfigure}	
	\caption{{The construction of orthogonal basis for synthetic functions $f_\mu(x)$. In the plots, u.p. RGS and m.p. RGS respectively refer to the unique precision RGS and the multi-precision RGS algorithms.}}
	\label{fig:Ex1_1}
\end{figure}

\Cref{fig:Ex1_1a} presents the evolution of the condition number of the computed Q factor at each iteration of GS process.  The evolution of (square root of) the condition number of $\bW_i$, $1\leq i \leq m$, is also depicted. We see that for $i \geq 150$, $\bW_i$ becomes numerically singular. For CGS and CGS2 methods, dramatic instabilities are observed at iterations $i \geq 50$ and $i \geq 150$, respectively.  MGS method exhibits more robustness than the other two standard variants of the GS process. With this method, the condition number of $\bQ_i$ remains close to $1$ up to iteration $i=130$, and then gradually degrades by more than an order of magnitude. The RGS algorithm executed in unique float32 format with $k=1500$ presents a similar stability as MGS. We see from~\cref{fig:Ex1_1a} that even though increasing of $k$ from $1500$ to $5000$ improves the quality of $\bTheta$ in terms of the $\varepsilon$-embedding property, the usage of $k=5000$ does not improve the stability of the RGS algorithm in unique float32 format but only worsens it. This can be explained by the increased rounding errors in computations of random projections and solutions of least-squares problems in Step 2. The multi-precision RGS algorithm, on the other hand, does not present this behavior.  It provides a Q factor with the condition number close to $1+\mathcal{O}(\varepsilon)$ and, particularly, an order of magnitude smaller than the condition number of the MGS Q factor.

The evolution of the approximation error $\|\bW_i - \bQ_i \bR_i\|/\|\bW_i\|$ is depicted in~\cref{fig:Ex1_1b}. We see that for CGS the error at first is close to the machine precision, but then it gradually degrades by two orders of magnitude. For CGS2 a dramatically large error is observed at iterations $i \geq 150$. For MGS and RGS algorithms the error remains close to the machine precision at all iterations.   

\Cref{fig:Ex1_1c} addresses a posteriori verification of the quality of the computed Q factor from its sketch.  We see that indeed $\cond(\bQ_i)$ can be well estimated by $\cond(\bS_i)$.  

Recall that the stability characterization of RGS algorithm in~\cref{stability} relies on the $\varepsilon$-embedding property of $\bTheta$. The minimal value $\omega$ of $\varepsilon$ for which $\bTheta$ satisfies the $\varepsilon$-embedding property for $\bQ_i$, at each iteration, is provided in \Cref{fig:Ex1_1d}. We also show the upper bound $\bar{\omega}$ for $\omega$ computed with~\cref{thm:cert01} from the sketches with no operations on high-dimensional vectors or matrices. In~\cref{thm:cert01}, the matrix $\bPhi$ was chosen to be of same size as $\bTheta$. Moreover, the parameter $\varepsilon^*$ was taken as $0.05$. It is observed that for both the unique precision and the multi-precision algorithms with $k \geq 5000$, the matrix $\bTheta$ satisfies the $\varepsilon$-embedding property with  (almost) $\varepsilon \leq 1/2$, which is the condition used in~\cref{stability} for deriving stability guarantees for RGS algorithm.   For $k = 1500$ at iterations $i\geq 70$, the value of $\omega$ becomes larger than $1/2$. Nevertheless, it remains small enough, which suggests a sufficient stability of RGS algorithm also for this value of $k$ and correlates well with the experiments (see~\cref{fig:Ex1_1a}). The estimator $\bar{\omega}$ of $\omega$ remains an upper bound of $\omega$ at all iterations and values of $k$, which implies robustness of~\cref{thm:cert01} for characterizing the $\varepsilon$-embedding property of $\bTheta$. An overestimation of $\omega$ by nearly a factor of $2$ is revealed at all iterations and values of $k$. Moreover, this behavior of $\bar{\omega}$ is observed also in other experiments. This suggests that, in practice, the value of $\bar{\omega}$ can be divided by a factor of $2$. 

\subsection{Orthogonalization of solution samples of a parametric PDE}
Next we consider a model order reduction problem from~\parencite[Section 6.1]{balabanov2019randomized2}. This problem describes a wave scattering with an object covered in an acoustic invisibility cloak. The cloak is multi-layered. The problem is governed by a parametric PDE, where the parameters are the properties of materials composing the last $10$ layers of the cloak, and the wave frequency.  By discretization with second-order finite elements, the parametric PDE is further transformed into a complex-valued system of equations of the form 
\begin{equation} \label{eq:exp2}
\bA_\mu \bu_\mu = \bb_\mu,
\end{equation}
where $\bA_\mu \in \mathbb{C}^{n \times n}$ and $\bb_\mu \in \mathbb{C}^{n}$ with $n \approx 400000$. The aim in~\parencite{balabanov2019randomized2} is to solve~\cref{eq:exp2} for parameters $\mu$ from the parameter set of interest $\mathcal{P}$. See~\parencite{balabanov2019randomized2} for more detailed description of the problem. 

Let us consider the construction of an orthogonal basis (so-called reduced basis) approximating the set $\{\bu_\mu: \mu \in \mathcal{P}\}$. For this, we drew from $\mathcal{P}$ $m=300$ uniform samples $\mu_1,\mu_2, \hdots, \mu_m$ and then performed a QR factorization of the matrix $$\bW = [\bu_{\mu_1}, \bu_{\mu_2}, \hdots, \bu_{\mu_m}]$$
with the following versions of the GS process: CGS, CGS2 and MGS computed in float32 format, the unique precision RGS in float32, and the multi-precision RGS using float32 for standard high-dimensional operations, while using float64 for other operations.\footnote{The extension of the theoretical analysis of RGS process from real numbers to complex numbers is straightforward.}

\begin{figure}[!h]
	\centering
	\begin{subfigure}{.35\textwidth}
		\centering  
		\includegraphics[width=\textwidth]{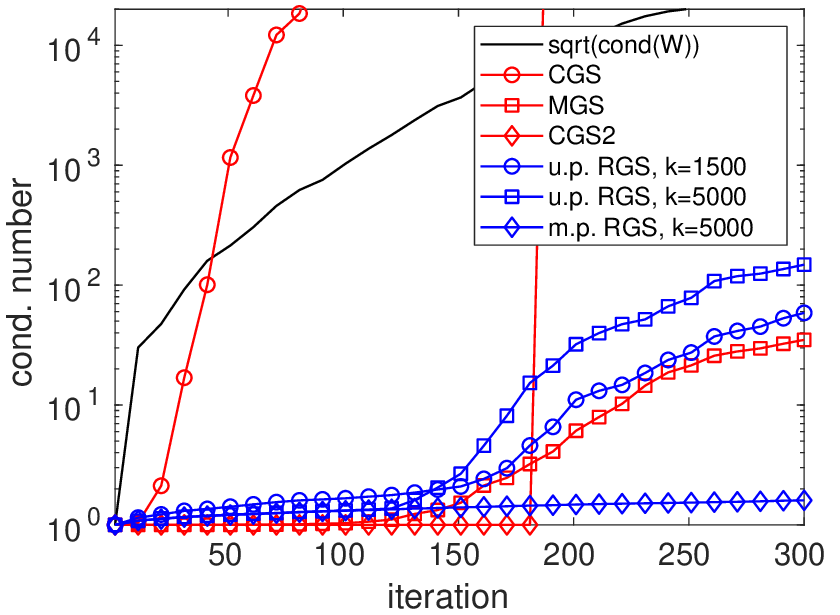}
		\caption{Cond. number of $\bQ_i$}
		\label{fig:Ex2_1a}
	\end{subfigure} \hspace{.1\textwidth}
	\begin{subfigure}{.35\textwidth}
		\centering
		\includegraphics[width=\textwidth]{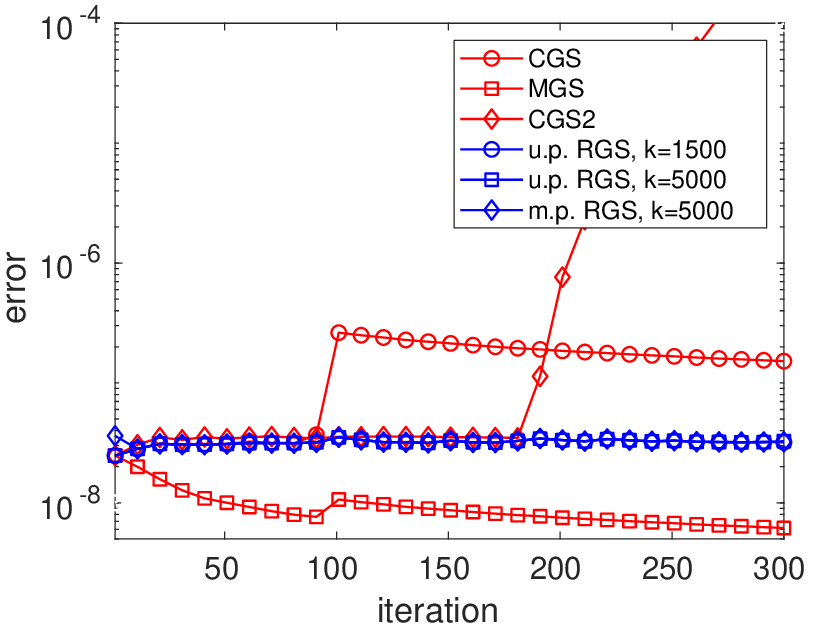}
		\caption{Error $\|\bW_i - \bQ_i \bR_i\|/\|\bW_i\|$}
		\label{fig:Ex2_1b}
	\end{subfigure} \\
	\begin{subfigure}{.35\textwidth}
		\centering  
		\includegraphics[width=\textwidth]{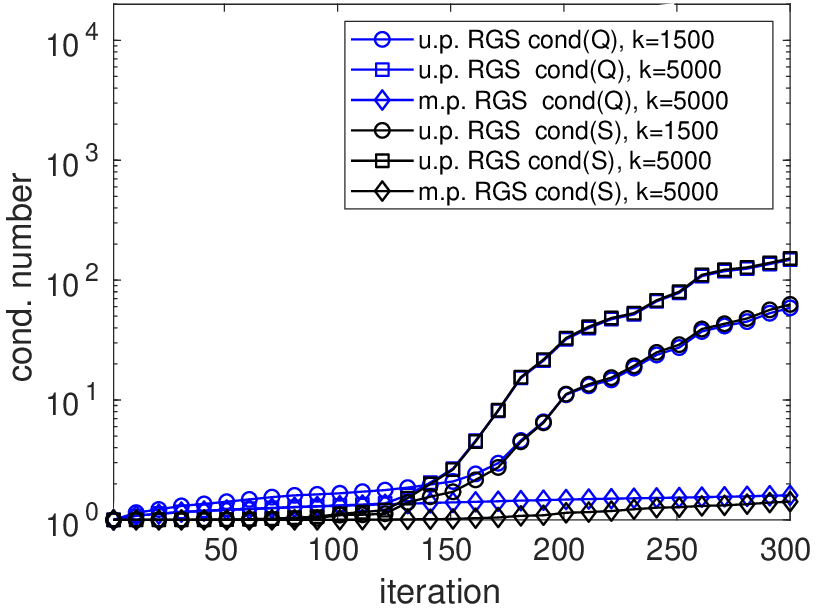}
		\caption{Cond. numbers of $\bQ_i$ and $\bS_i$}
		\label{fig:Ex2_1c}
	\end{subfigure} \hspace{.1\textwidth}
	\begin{subfigure}{.35\textwidth}
		\centering
		\includegraphics[width=\textwidth]{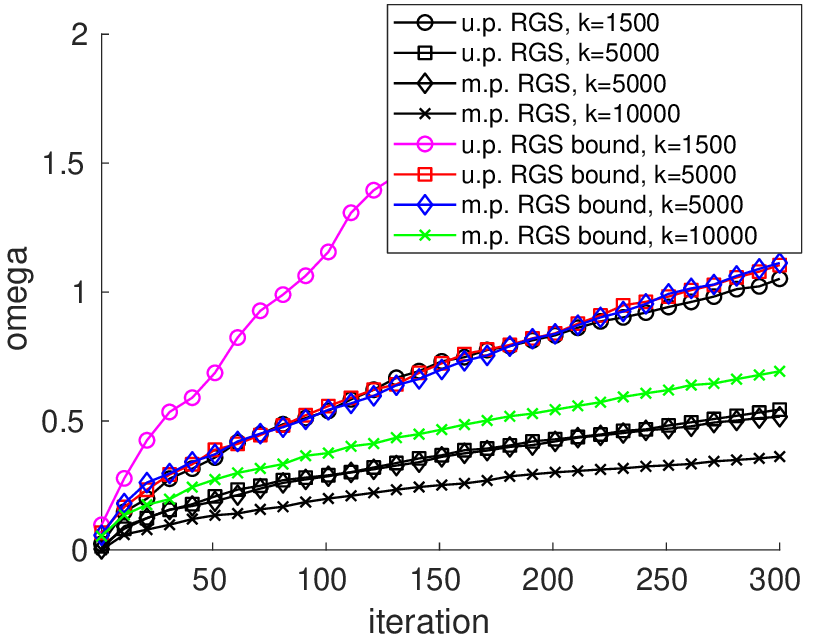}
		\caption{Coeff. $\omega$ and its upper bound}
		\label{fig:Ex2_1d}
	\end{subfigure}	
	\caption{{The construction of orthogonal basis for the solution set of a parametric PDE.  In the plots, u.p. RGS and m.p. RGS respectively refer to the unique precision RGS and the multi-precision RGS algorithms.}}
	\label{fig:Ex2_1}
\end{figure}

The condition number of the factor $\bQ_i$ and the approximation error $\|\bW_i - \bQ_i \bR_i\|/\|\bW_i\|$ obtained at each iteration of the algorithms are depicted in~\cref{fig:Ex2_1a,fig:Ex2_1b}, respectively. Furthermore, for randomized algorithms, in~\cref{fig:Ex2_1c} we provide a comparison of $\cond(\bQ_i)$ and $\cond(\bS_i)$. In~\cref{fig:Ex2_1d} we present the characterization of the $\varepsilon$-embedding property of $\bTheta$ for $\bQ_i$ given by the value $\omega$ and its upper bound computed with~\cref{thm:cert01}. In~\cref{thm:cert01} we chose $\bPhi$ to be of same size as $\bTheta$ with the parameter $\varepsilon^* = 0.05$.

A very similar picture is observed as in the previous numerical example.  More specifically, dramatic instabilities are revealed at iterations $i \geq 50$ for CGS and $i \geq 190$ for CGS2. The MGS and the unique precision RGS with $k=1500$ show a similar stability, which is better than the one of CGS and CGS2.  For these algorithms $\cond(\bQ_i)$ remains close to $1$ at iterations 
$i \leq 150$, but degrades by more than an order of magnitude at latter iterations.  
The multi-precision randomized GS algorithm at all iterations provides a Q factor with condition number close to $1$. The approximation error $\|\bW_i - \bQ_i \bR_i\|/\|\bW_i\|$ again is close to the machine precision for MGS and RGS algorithms, while it is larger for CGS and CGS2. The condition number of the sketched Q factor, $\bS_i$, is verified to be a good estimator of the condition number of $\bQ_i$. For both the unique as well as the multi-precision RGS algorithms with $k \geq 5000$, the condition $\omega \leq 1/2$ used in the stability analysis is (nearly) satisfied. For $k=1500$, the value of $\omega$ is larger than $1/2$ at iterations $i \geq 70$,  though it remains sufficiently small suggesting the stability of RGS algorithm also for this sketch size. Again, we  reveal an overestimation of $\bar{\omega}$ as an upper bound of $\omega$, by nearly a factor of 2. 

\subsection{Solution of a linear system with GMRES} \label{expgmres}
In this numerical experiment the RGS algorithm is tested in the context of GMRES method for the solution of the linear system of equations: 
\begin{equation}  \label{eq:ex3sys}
\bA_{f} \bx_{f} = \bb,
\end{equation}
where the matrix $\bA_f$ is taken as the ``SiO2'' matrix of dimension $n = 155331$ from the {SuiteSparse  matrix  collection}. The right-hand-side vector $\bb$ is taken as $\bb = \bA \by/\|\bA \by \|$, where  $\by=  [1,1, \hdots ,1 ]^\mathrm{T}$. Furthermore, the system~\cref{eq:ex3sys} is preconditioned from the right by the incomplete LU factorization $\bP_f$ of $\bA$ with zero level of fill in. With this preconditioner the final system of equations has the following form$$  \bA \bx = \bb, $$ 
where $\bA = \bA_f \bP_f$ and $\bx_{f} = \bP_f \bx $. This system is considered for the solution with GMRES method based on different versions of GS process. Here we test only CGS, CGS2, MGS and the unique precision RGS. The multi-precision RGS is not considered since the unique precision RGS already provides a nearly optimal solution. 

In all experiments, the products with matrix $\bA$ are computed in float64 format. The solutions of the Hessenberg least-squares problems~\cref{eq:hls} are computed with Givens rotations also in float64 format. All other operations (i.e., the GS iterations) are performed in float32 format. 

\begin{figure}[!h]
	\centering
	\begin{subfigure}{.35\textwidth}
		\centering  
		\includegraphics[width=\textwidth]{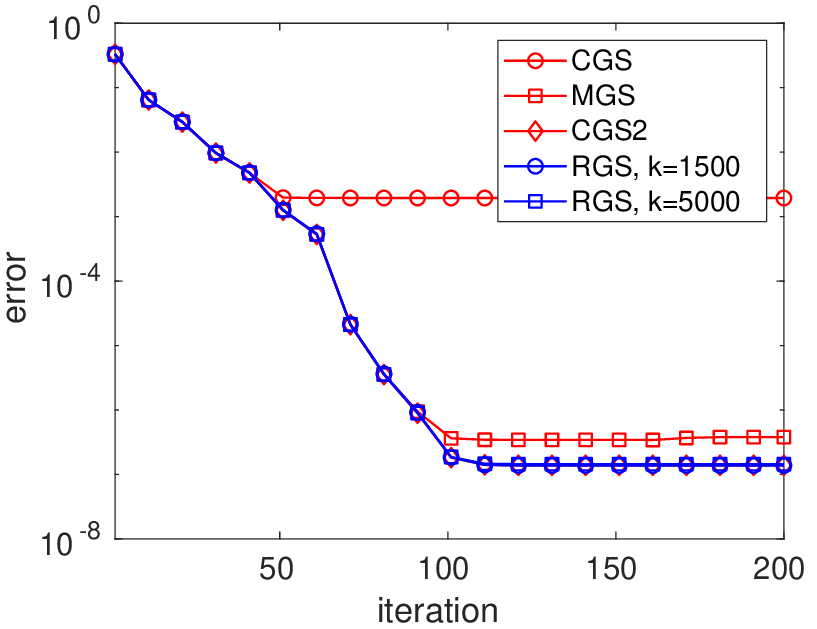}
		\caption{Residual error $\|\bA \bhx_{i-1} - \bb \|$.}
		\label{fig:Ex3_1}
	\end{subfigure} \hspace{.1\textwidth}
	\begin{subfigure}{.35\textwidth}
		\centering
		\includegraphics[width=\textwidth]{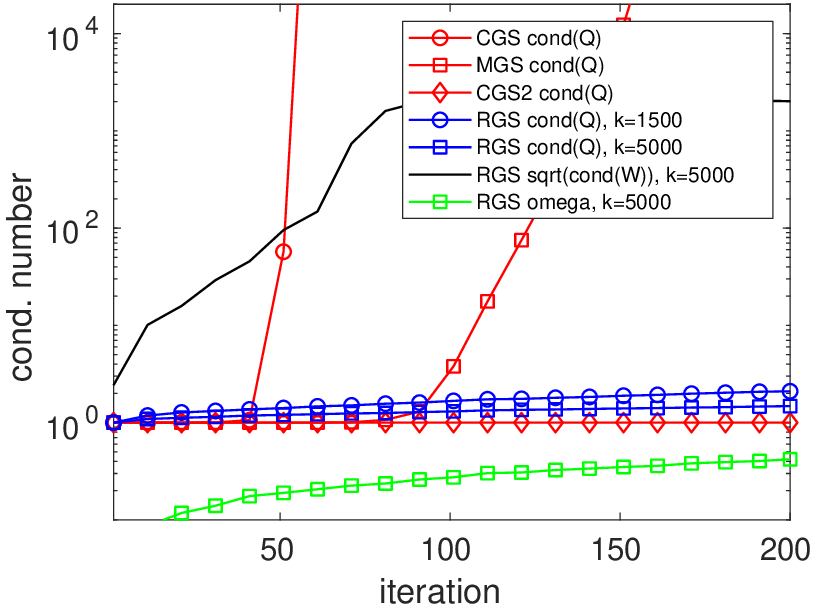}
		\caption{Cond. number of $\bQ_i$}
		\label{fig:Ex3_2}
	\end{subfigure}
	\caption{Solution of a linear system with GMRES.}
	\label{fig:Ex3}
\end{figure}

The convergence of the residual error is depicted in~\cref{fig:Ex3_1}. The condition number of the  Q factor (characterizing the orthogonality of the computed Krylov basis) at each iteration $i$ is provided in~\cref{fig:Ex3_2}. In~\cref{fig:Ex3_2} we also provide the condition number of $\bW_i = [\bA\bQ_{i-1},~\bb]$ and the value of $\omega$ representing the $\varepsilon$-embedding property of $\bTheta$ for $\bQ_i$ in the RGS algorithm with $k =5000$. 

At iterations $i \geq 50$, we reveal a dramatic instability of the CGS algorithm resulting in the early stagnation of the residual error. The other versions of the GS process present a better stability. The CGS2 and RGS with all sketch sizes at all iterations yield almost orthogonal Q factor. For these algorithms the error has converged to machine precision. The MGS algorithm does not provide a well-conditioned Q factor at $i \geq 110$ iterations. Nevertheless, it yields the convergence of the residual error up to machine precision similar to CGS2 and RGS.  

Finally, we see that for RGS with $k=5000$, $\bTheta$ is verified to be an $\varepsilon$-embedding for $\bQ_i$ with $\varepsilon \leq 1/2$. This implies applicability of the stability analysis from~\cref{stability}. 

\section{Conclusion} \label{concl}

In this article we proposed a novel randomized  Gram-Schmidt  process for efficient orthogonalization of a set of high-dimensional vectors. This process can be incorporated into GMRES method or Arnoldi iteration for solving large systems of equations or eigenvalue problems.  Our methodology can be adapted to practically any computational architecture. 

The randomized  GS process was introduced under multi-precision arithmetic model that also accounts for the classical unique precision model. We proposed to perform expensive high-dimensional operations in low precision, while computing the inexpensive random projections and low-dimensional operations in high precision. The numerical stability of the algorithms was shown for the low-precision unit roundoff independent of the dimension of the problem. This feature can have a major importance when solving extreme-scale problems.  

{The great potential of the methodology was realized with three numerical examples. In all the experiments, the multi-precision RGS algorithm provided a Q factor (i.e., the orthogonalized matrix) with condition number close to $1$. It remained stable even in extreme cases, such as orthogonalization of a numerically singular matrix, where the standard CGS and CGS2 methods failed.  This is in addition to the fact that the RGS algorithm can require nearly {half as many} flops and passes over the data than CGS.} 

\section{Acknowledgments}
This project has received funding from the European Research Council (ERC) under the European Union's Horizon 2020 research and innovation program (grant agreement No 810367).

\newpage

\printbibliography

\newpage
\section*{Supplementary materials}
\renewcommand*{\thesection}{SM}
\begin{refsection}
	
	Here we provide the proofs of some propositions and theorems from the paper.
	
	\begin{proof} [Proof of~\Cref{thm:thetadeltax0}]	
		{Our proof is based on that of \parencite[Theorem 3.1]{krahmer2011new}.} 	We first notice that it is sufficient to prove the theorem for a normalized $\bgamma$. Moreover, we may permute the entries of $\bgamma$ and $\bphi$ in any convenient way, since the RIP property is invariant under permutation of columns of $\bTheta$. Consequently, without loss of generality we further assume that the entries of $\bgamma$ are ordered in a non-decreasing way and $\|\bgamma\|=1$. 
		
		For a vector $\by \in \mathbb{R}^n$ and index $ I \in \{1,2, \hdots,  \lceil \frac{n}{d} \rceil \}$  we denote by $\by_{(I)} \in \mathbb{R}^n$  a vector with entries 
		$$[\by_{(I)}]_i = \left \{
		\begin{array}{rll}
		&[\by]_i,~~~~  \textup{if }  d(I-1)+1 \leq i \leq \min(dI,n) \\
		&0,~~~\textup{otherwise.} 
		\end{array}
		\right.$$ 
		The RIP property of $\bTheta$ and the parallelogram identity then imply that
		\begin{equation}\label{eq:thetadeltax1}
		\langle \bTheta \by_{(I)},  \bTheta \bz_{(J)} \rangle \leq (\varepsilon/4) \|\by_{(I)}\| \|\bz_{(J)}\|,
		\end{equation}
		holds for any vectors $\by,\bz \in \mathbb{R}^n$ and $I,J = 1,2, \hdots, \lceil \frac{n}{d} \rceil$, $I \neq J$. {This classical result can be found for instance in~\parencite{krahmer2011new}}.
		
		We have $$\bphi = \mathrm{diag}(\bxi) \bgamma = \mathrm{diag}(\bgamma) \bxi,$$ 
		where $\bxi = (\xi_1, \xi_2, \hdots, \xi_n)$ is a vector with entries that are independent random variables in~$[-1,1]$. Furthermore, 
		\begin{equation} \label{eq:Thetaphisum}
		\begin{split}
		\| \bTheta \bphi \|^2  
		&= \| \bTheta \mathrm{diag}(\bgamma) \bxi \|^2 \\
		&=  \sum^{\lceil \frac{n}{d} \rceil}_{I=1} \| \bTheta \mathrm{diag}(\bgamma_{(I)}) \bxi \|^2 + 2  \sum^{\lceil \frac{n}{d} \rceil}_{I=2} \langle \bTheta \mathrm{diag}(\bgamma_{(I)})\bxi , \bTheta \mathrm{diag}(\bgamma_{(1)}) \bxi \rangle \\ 
		&+ \sum^{\lceil \frac{n}{d} \rceil}_{I,J = 2,  I \neq J} \langle  \bTheta \mathrm{diag}(\bgamma_{(I)})\bxi , \bTheta \mathrm{diag}(\bgamma_{(J)}) \bxi \rangle \\ 
		&=  \sum^{\lceil \frac{n}{d} \rceil}_{I=1} \| \bTheta \mathrm{diag}(\bgamma_{(I)}) \bxi \|^2 + 2 \langle \bv, \bxi \rangle + \langle \bC \bxi, \bxi  \rangle,
		\end{split}
		\end{equation}
		where $\bv := \sum^{\lceil \frac{n}{d} \rceil}_{I=2} \mathrm{diag}(\bgamma_{(I)}) \bTheta^\mathrm{T} \bTheta \mathrm{diag}(\bgamma_{(1)}) \bxi_{(1)}$ and $\bC : = \sum^{\lceil \frac{n}{d} \rceil}_{I,J = 2,  I \neq J} \mathrm{diag}(\bgamma_{(I)}) \bTheta^\mathrm{T} \bTheta \mathrm{diag}(\bgamma_{(J)})$. 
		
		Each term in~\cref{eq:Thetaphisum} can be estimated separately: 
		\begin{itemize}
			\item Since $\bTheta$ satisfies $(\varepsilon/4, 2d)$-RIP, we have
			$$ (1-\varepsilon/4) \|\mathrm{diag}(\bgamma_{(I)}) \bxi \|^2  \leq \| \bTheta \mathrm{diag}(\bgamma_{(I)})  \bxi \|^2 \leq (1+\varepsilon/4) \|\mathrm{diag}(\bgamma_{(I)})  \bxi \|^2,~1 \leq I \leq \lceil \frac{n}{d} \rceil,$$
			which implies that 
			$$ (1-\varepsilon/4) \| \bphi \|^2 \leq   \sum^{\lceil \frac{n}{d} \rceil}_{I=1} \| \bTheta \mathrm{diag}(\bgamma_{(I)}) \bxi \|^2 \leq  (1+\varepsilon/4) \| \bphi \|^2.$$ 
			\item Notice that 
			\begin{equation} \label{eq:thetadeltax2}
			\| \bgamma_{(I)} \|_\infty \leq \frac{1}{\sqrt{d}}\| \bgamma_{(I-1)} \|_2, \text{ for}~I\geq 2. 
			\end{equation}		
			Consequently, by~\cref{eq:thetadeltax1,eq:thetadeltax2} and the fact that $\|\mathrm{diag}(\bgamma_{(1)}) \bxi_{(1)} \| \leq \| \bgamma \| = 1,$ we have
			\small 
			\begin{equation*}
			\begin{split}
			\|\bv\| 
			&= \| \sum^{\lceil \frac{n}{d} \rceil}_{I=2}   \mathrm{diag}(\bgamma_{(I)}) \bTheta^\mathrm{T} \bTheta \mathrm{diag}(\bgamma_{(1)}) \bxi_{(1)}\| = \sup_{\|\by\|=1} \sum^{\lceil \frac{n}{d} \rceil}_{I=2} \by^\mathrm{T}  \mathrm{diag}(\bgamma_{(I)}) \bTheta^\mathrm{T} \bTheta \mathrm{diag}(\bgamma_{(1)}) \bxi_{(1)} \\
			&  \leq    \sup_{\|\by\|=1} \sum^{\lceil \frac{n}{d} \rceil}_{I=2} (\varepsilon/4) \| \mathrm{diag}(\bgamma_{(I)}) \by \|  \|\mathrm{diag}(\bgamma_{(1)}) \bxi_{(1)} \| \leq \sup_{\|\by\|=1} \sum^{\lceil \frac{n}{d} \rceil}_{I=2} (\varepsilon/4) \| \bgamma_{(I)}\|_\infty \| \by_{(I)} \|  \\ 
			&  \leq  (\varepsilon/4)  \sup_{\|\by\|=1} \sum^{\lceil \frac{n}{d} \rceil}_{I=2} \frac{1}{\sqrt{d}} \| \bgamma_{(I-1)}\| \|\by_{(I)} \|    \leq \frac{\varepsilon}{4\sqrt{d}}   \sup_{\|\by\|=1} \sum^{\lceil \frac{n}{d} \rceil}_{I=2} (\| \bgamma_{(I-1)}\|^2 + \| \by_{(I)} \|^2)/2  \\
			&\leq \frac{\varepsilon}{4\sqrt{d}}.
			\end{split}
			\end{equation*}
			\normalsize
			
			Define $\bxi':= \sum^{\lceil \frac{n}{d} \rceil}_{I=2} \bxi_{(I)}$. Notice that $\bxi'$ and $\bv$ are randomly independent and $\langle \bv, \bxi \rangle= \langle \bv, \bxi' \rangle $.
			The Hoeffding's inequality states that for any $t >0$,		
			$$ \mathbb{P} \left ( |\langle \bv, \bxi' \rangle| \geq t \right )  \leq 2 \exp \left (-\frac{t^2}{2\| \bv\|^2} \right ), $$ 
			from which we conclude that
			$$ \mathbb{P} \left ( |\langle \bv, \bxi \rangle| \geq \varepsilon/4 \right )  \leq 2 \exp (-d/2). $$ 
			\item Let $\bone$ denote a vector with entries $[\bone]_i =1$, $1 \leq i \leq n$. The estimation of the last term in~\cref{eq:Thetaphisum}, shall rely on the following bounds for $\|\bC\|$ and $\|\bC\|_\Frob$ derived by using relations~\cref{eq:thetadeltax1,eq:thetadeltax2}:	
			\begin{equation*}
			\begin{split}		
			\| \bC \| 
			&=  \sup_{\|\by\|=1} \langle \by, \bC \by \rangle \\
			&=\sup_{\|\by\|=1} \sum^{\lceil \frac{n}{d} \rceil}_{I,J = 2,  I \neq J} \langle \bTheta \mathrm{diag}(\bgamma_{(I)}) \by_{(I)},  \bTheta \mathrm{diag}(\bgamma_{(J)}) \by_{(J)} \rangle \\
			&\leq \sup_{\|\by\|=1} \sum^{\lceil \frac{n}{d} \rceil}_{I,J = 2,  I \neq J} (\varepsilon/4) \| \mathrm{diag}(\bgamma_{(I)}) \by_{(I)} \|  \|\mathrm{diag}(\bgamma_{(J)}) \by_{(J)}\| \\
			& \leq (\varepsilon/4)	 \sup_{\|\by\|=1} \sum^{\lceil \frac{n}{d} \rceil}_{I,J = 2,  I \neq J}  \| \by_{(I)} \| \| \by_{(J)} \| \|\bgamma_{(I)}\|_{\infty} \|\bgamma_{(J)}\|_{\infty} \\
			& \leq (\varepsilon/4)	 \sup_{\|\by\|=1} \sum^{\lceil \frac{n}{d} \rceil}_{I,J = 2,  I \neq J}  \| \by_{(I)} \| \| \by_{(J)} \| \frac{1}{\sqrt{d}}\| \bgamma_{(I-1)} \| \frac{1}{\sqrt{d}}\| \bgamma_{(J-1)} \| \\
			& \leq \frac{\varepsilon}{4d}	 \sup_{\|\by\|=1} \sup_{\|\by\|=1} \sum^{\lceil \frac{n}{d} \rceil}_{I,J = 2,  I \neq J} \left ( \| \by_{(I)} \|^2 +  \| \bgamma_{(I-1)} \|^2 \right ) \left ( \| \by_{(J)} \|^2 +  \| \bgamma_{(J-1)} \|^2 \right )/4 \\
			& \leq \frac{\varepsilon}{4d},	
			\end{split}
			\end{equation*}
			and, 
			\begin{equation*}
			\begin{split}		
			\| \bC \|^2_\Frob  
			&=\sum^{\lceil \frac{n}{d} \rceil}_{I,J = 2,  I \neq J} \| \mathrm{diag}(\bgamma_{(I)}) \bTheta^\mathrm{T} \bTheta \mathrm{diag}(\bgamma_{(J)}) \|^2_\Frob \\
			& \leq \sum^{\lceil \frac{n}{d} \rceil}_{I,J = 2,  I \neq J} \| \mathrm{diag}(\bgamma_{(I)})\|^2 \| \mathrm{diag}(\bone_{(I)}) \bTheta^\mathrm{T} \bTheta \mathrm{diag}(\bone_{(J)})\|^2 \| \mathrm{diag}(\bgamma_{(J)}) \|^2_\Frob \\
			& \leq \sum^{\lceil \frac{n}{d} \rceil}_{I,J = 2,  I \neq J} \| \bgamma_{(I)}\|^2_\infty \left ( \max_{\|\bx\| =1, \|\by\| =1} \langle \bTheta \by_{(I)},\bTheta \bx_{(J)} \rangle \right)^2 \|\bgamma_{(J)} \|^2 \\
			& \leq (\varepsilon/4)^2 \sum^{\lceil \frac{n}{d} \rceil}_{I,J = 2,  I \neq J} 
			\|\bgamma_{(I)} \|_\infty^2 \|\bgamma_{(J)}\|^2 \\
			& \leq  (\varepsilon/4)^2 \sum^{\lceil \frac{n}{d} \rceil}_{I=2} \frac{1}{d}
			\|\bgamma_{(I-1)} \|^2 \\
			&\leq 
			\left (\frac{\varepsilon}{4\sqrt{d}} \right )^2.	
			\end{split}
			\end{equation*}
			
			Now, the Hanson-Wright inequality (see~\parencite[Theorem 6.2.1]{vershynin2018high}) states that for any $t >0$,		
			$$ \mathbb{P} \left ( |\langle \bxi, \bC \bxi \rangle| \geq t \right )  \leq 2 \exp \left (- c \min (\frac{t}{K^2 \| \bC \|}, \frac{t^2}{ K^4\| \bC \|_\Frob^2} ) \right ), $$ 
			where $K = \max_i \| \xi_i \|_{\psi_2}$ and $c$ is a universal constant. Since  the random variables $\xi_i$ lie in the interval $[-1,1]$, hence their sub-Gaussian norms satisfy $\| \xi_i \|_{\psi_2} \leq 1/\sqrt{\log(2)}$. From this fact and the derived earlier bounds for $\| \bC \|$ and $\| \bC \|_\Frob$ we deduce that 
			$$ \mathbb{P} \left ( |\langle \bxi, \bC \bxi \rangle| \geq \varepsilon/4 \right )  \leq 2 \exp \left (- c \min (\log(2) d ,  \log^2(2) d  ) \right ), $$ 
			which results in
			$$ \mathbb{P} \left ( |\langle \bxi, \bC \bxi \rangle| \geq \varepsilon/4 \right )  \leq 2 \exp (- 0.48 c d). $$ 			
		\end{itemize}
		
		We conclude that with probability at least $1-(2 \exp (-d/2)+2 \exp (- 0.48 c d))$:
		\begin{equation*}
		\begin{split}
		| \|\bTheta \bphi  \|^2 - \|\bphi\|^2 | &= 	
		\left | \sum^R_{J=1} \| \bTheta \mathrm{diag}(\bgamma_{(J)}) \bxi \|^2 + 2 \langle \bv, \bxi \rangle + \langle \bC \bxi, \bxi  \rangle - \|\bphi\|^2 \right |  \\
		&\leq \left| \sum^R_{J=1} \| \bTheta \mathrm{diag}(\bgamma_{(J)}) \bxi \|^2 - \|\bphi\|^2 \right |+ 2 |\langle \bv, \bxi \rangle| + |\langle \bC \bxi, \bxi  \rangle| \\
		& \leq \varepsilon/4+ 2 \varepsilon/4 + \varepsilon/4 = \varepsilon.
		\end{split}
		\end{equation*}
		The proof is finished by noting that $\delta < (2 \exp (-d/2)+2 \exp (- 0.48 c d))$.
	\end{proof}

	\begin{proof}[Proof of~\Cref{thm:maintheorem1}]
		By standard worse-case scenario rounding analysis and the assumptions~\cref{eq:Thetanorm,eq:numepsembedding3}, we have for $1 \leq i \leq m$:
		\begin{align*}
		\| \bTheta \bhq_i - \bhs_i\| &\leq \| \bTheta(\bhq_i - \bhq_i'/\hat{r}_{i,i}) \| + \| \bTheta \bhq_i' - \bhs_i' \|/|\hat{r}_{i,i}| + \| \bhs_i - \bhs_i'/\hat{r}_{i,i} \|\\
		&\leq  u_{fine}  D  \| \bhq_i'/\hat{r}_{i,i} \|  + 1.02 n u_{fine} \| \bTheta \|_\Frob \| \bhq_i'/\hat{r}_{i,i} \| + u_{fine} \| \bhs_i'/\hat{r}_{i,i} \| \\
		&\leq 1.01( 10^{-4} u_{crs} m^{-1/2} \| \bhq_i\|  + 0.0102 u_{crs}  \sqrt{3/2} m^{-1/2} \| \bhq_i\| +10^{-5} u_{crs}  m^{-1/2}) \\
		& \leq 0.013 u_{crs} m^{-1/2} \| \bhq_i\| + 10^{-4} u_{crs} m^{-1/2}.
		\end{align*}
		Moreover, this implies that
		\small 
		\begin{equation} \label{eq:main10}
		\begin{split}
		&\| \bTheta \bhQ_m - \bhS_m\|^2_\Frob  = \sum^m_{i=1} \| \bTheta \bhq_i - \bhs_i\|^2 \leq \sum^m_{i=1} (0.013 u_{crs} m^{-1/2} \| \bhq_i\| + 10^{-4} u_{crs} m^{-1/2})^2 \\&			
		\leq \sum^m_{i=1}2( 0.013^2 u^2_{crs} m^{-1} \| \bhq_i\|^2 + 10^{-8}   u^2_{crs} m^{-1}) = 2( 0.013^2 u^2_{crs} m^{-1} \| \bhQ_m\|_\Frob^2 + 10^{-8}   u^2_{crs}) \\
		&\leq {(}0.02 u_{crs} m^{-1/2}\| \bhQ_m\|_\Frob + 0.01 u_{crs}{)^2}  \leq u_{crs}{^2} (0.02  \| \bhQ_m\| + 0.01){^2}:=F_1{^2}.
		\end{split}
		\end{equation}
		\normalsize
		
		In its turn, this implies that			
		$$ \sigma_{min}(\bhS) - F_1 \leq \sigma_{min}(\bTheta\bhQ) \leq \sigma_{max}(\bTheta\bhQ) \leq \sigma_{max}(\bhS) + F_1.$$
		This result combined with the fact that
		$$ 1-\Delta_m \leq \sqrt{1-\Delta_m} \leq  \sigma_{min}(\bhS)  \leq \sigma_{max}(\bhS) \leq \sqrt{1+\Delta_m} \leq 1+\Delta_m ,$$
		and~\cref{thm:epscond} imply that $\|\bhQ \| \leq 1.575$ and that
		$$ (1+\varepsilon)^{-1/2}(1-\Delta_m- F_1) \leq \sigma_{min}(\bhQ)  \leq \sigma_{max}(\bhQ) \leq (1-\varepsilon)^{-1/2}(1+\Delta_m+ F_1).$$
		
		Let us now prove the second statement of the theorem. Assume that $m\geq 2$ (otherwise the proof is obvious). We first notice that
		\begin{equation}	\label{eq:main11}
		\begin{split}
		\| \bhR \|_\Frob 
		&\leq \sigma_{min}(\bhS)^{-1} \| \bhS \bhR\|_\Frob \leq 1.12 \| \bhS \bhR\|_\Frob \leq 1.12 (\|\bhP \|_\Frob +  \|\bhP - \bhS \bhR\|_\Frob ) \\ 
		&\leq  1.12 (1+\tilde{\Delta}_m)\|\bhP\|_\Frob \leq 1.4 (1+\tilde{\Delta}_m)\|\bhW\|_\Frob  \leq 1.6 \|\bhW\|_\Frob. 
		\end{split}
		\end{equation}
		In~\cref{eq:main11} we used the fact that $\| \bhp_i \| \leq 1.02 \| \bTheta \bhw_i \| \leq 1.25 \| \bhw_i\|$ so that $ \|\bhP\|_\Frob \leq 1.25 \| \bhW\|_\Frob$.
		Furthermore, by using the relation~\cref{eq:numepsembedding3stand}, we obtain
		\small
		\begin{equation*}	
		\begin{split}
		&\|\bhw_i - \bhQ [\bhR]_{(1:m,i)}\| 
		\leq \|\bhw_i - \bhq_i' - \bhQ_{i-1} [\bhR]_{(1:i-1,i)}\| + \| \bhq_i' - \bhq_i r_{i,i} \|\\
		&\leq \|\bDelta \bhq'_i\| +u_{fine}  \|\bhq'_i \| 
		\leq 1.02 u_{crs}(\| \bhw_i \| +  i \|\bhQ_{i-1}\|_\Frob \|\bhr_{i} \|) +1.1 u_{fine}  \|\bhq_i \| \|\bhr_i\|\\
		& \leq 1.02 u_{crs}\| \bhw_i \| + 1.61  u_{crs} i^{3/2} \|\bhr_i\|.
		\end{split}
		\end{equation*}
		\normalsize
		Consequently, 
		\small
		\begin{align*}
		&\| \bhW - \bhQ \bhR \|^2_\Frob 
		= \sum^m_{i=1} \|\bhw_i - \bhQ [\bhR]_{(1:m,i)}\|^2 
		\leq \sum^m_{i=1} u^2_{crs} (1.02 \|\bhw_i \| +1.61 m^{3/2}\|\bhr_i\|)^2 \\
		& \leq \sum^m_{i=1} u^2_{crs} 2(1.02^2 \|\bhw_i \|^2 +1.61^2 m^{3} \|\bhr_i\|^2)  
		\leq u^2_{crs} 2(1.02^2 \|\bhW \|^2_\Frob +1.61^2 m^{3} \|\bhR\|_\Frob^2) \\
		&\leq u^2_{crs} 2(1.02^2 +2.58^2 m^{3}) \|\bhW\|_\Frob^2
		\leq \left ( 3.7 u_{crs} m^{3/2} \|\bhW\|_\Frob \right)^2,
		\end{align*}	
		\normalsize
		which finishes the proof.
	\end{proof}

	\begin{proof}[Proof of~\cref{thm:maintheorem2}]
		The proof is done by induction on $m$. The inequalities~\cref{eq:maintheorem21,eq:maintheorem22} of the theorem are obvious for $m=1$. Now, suppose that they are true for $m=i-1$. The goal is to derive the inequalities for $m=i \geq 2$. 
		
		{From~\parencite[Theorems 8.5, 19.10 and 20.3]{higham2002accuracy} and their proofs,} the least squares solution in Step 2 satisfies the following backward-stability property:
		\begin{equation} \label{eq:bls}
		[\bhR]_{(1:i-1,i)}  = \arg \min_{\by} \| (\bhS_{i-1}+ \bDelta\bS_{i-1})\by - (\bhp_i + \bDelta\bp_i) \|,~\text{ with}
		\end{equation}
		$$   \| \bDelta\bS_{i-1} \|_\Frob \leq 0.01 u_{crs}   \| \bhS_{i-1} \|,~~   \| \bDelta\bp_i \| \leq 0.01 u_{crs}   \|\bhp_i \|,$$ 
		
		Clearly,
		\begin{equation} \label{eq:mainproof1}
		{\|\bhS_i\|_\Frob = (\sum^{i}_{j=1} \|\bhs_j\|^2)^{1/2} \leq 1.01 \sqrt{i} \text{~ and ~} \| \bhp_i \| \leq 1.02 \|\bTheta \bhw_i \| \leq {1.25} \|\bhw_i\|.}
		\end{equation}
		Moreover, we have for $1\leq j \leq i-1$ {(which can be shown to hold similarly to~\cref{eq:main10} in the proof of \cref{thm:maintheorem1} taking $m=j$)},
		\begin{equation} \label{eq:ThetaQ}
		\|\bhS_j - \bTheta \bhQ_j\|_\Frob \leq  u_{crs}   (0.02\|\bhQ_j\|+0.01).
		\end{equation}	
		{From the fact that $\sigma_{min}(\bhS_{i-1}) \geq \sqrt{1-\Delta_{i-1}} \geq 0.989 $, and~\cref{eq:bls,eq:mainproof1}, notice that} 
		\begin{equation} \label{eq:mainproof2}
		\sigma_{min}(\bhS_{i-1}+ \bDelta\bS_{i-1}) \geq  \sigma_{min}(\bhS_{i-1})- \| \bDelta \bhS_{i-1}\| \geq 0.98. 
		\end{equation}
		Consequently, we have		
		\begin{equation} \label{eq:mainproof3}
		\| [\bhR]_{(1:i-1,i)}\| \leq \| \bhp_i + \bDelta \bp_i \|/\sigma_{min}(\bhS_{i-1}+ \bDelta\bS_{i-1}) \leq 1.4 \| \bhw_i \|.
		\end{equation} 	
		{This fact combined with the fact (which holds by the induction hypothesis: $\Delta_{i-1} \leq 0.02$ and~\cref{thm:maintheorem1}) that} 
		\begin{equation} \label{eq:mainproof4}
		\|\bhQ_{i-1}\| \leq 1.5 \text{ and } \|\bhQ_{i-1}\|_\Frob \leq 1.5 i^{1/2}, 
		\end{equation}
		and~\cref{eq:numepsembedding3stand,eq:mainproof3,eq:mainproof4}, leads to the following result: 
		\small	
		\begin{gather}  
		\|\bhq'_i-(\bhw_i - \bhQ_{i-1} [\bhR]_{(1:i-1,i)}) \| \leq 1.02 u_{crs}  \| |\bhw_i|+ i | 	\bhQ_{i-1}| |[\bhR]_{(1:i-1,i)} | \| \leq 3.2 u_{crs} i^{3/2} \|\bhw_i \|, \label{eq:qprime} \\ 	
		\text{ and ~~~~~}\| \bhq'_i \| \leq \|\bhw_i\|+ \|\bhQ_{i-1}\| \|[\bhR]_{(1:i-1,i)}\| + \|\bDelta \bhq'_i\| \leq 3.2 i^{1/2} \|\bhw_i\|. \label{eq:qprime2}
		\end{gather}	
		Denote $\bhq'_i-(\bhw_i - \bhQ_{i-1} [\bhR]_{(1:i-1,i)})$ by $\bDelta\bq'_i$.
		
		In addition to~\cref{eq:qprime,eq:qprime2}, we also have by~\cref{eq:mainproof3,eq:mainproof4,eq:numepsembedding3stand,eq:numepsembedding3} that
		\begin{gather} 
		\bTheta \bhq'_i =  \bTheta \bhw_i - \bTheta \bhQ_{i-1} [\bhR]_{(1:i-1,i)} + \bTheta \bDelta\bq'_i, \text{with} \label{eq:mainproof5} 	\\
		\|\bTheta \bDelta\bq'_i \|\leq  1.02 u_{crs} D \| |\bhw_i| + i | \bhQ_{i-1}| |[\bhR]_{(1:i-1,i)}| \| \leq  4 u_{crs}   i^{3/2}  \| \bhw_i \|. \nonumber
		\end{gather}		 
		By combining~\cref{eq:mainproof5} and \cref{eq:ThetaQ} along with standard rounding analysis and \cref{eq:qprime2,eq:ThetaQ,eq:Thetanorm,eq:mainproof3,eq:mainproof4,eq:mainproof5}, it follows that	
		\begin{equation} \label{eq:sprime}
		\begin{split}
		&\bhs'_i = \bhp_i - \bhS_{i-1} [\bhR]_{(1:i-1,i)} + \bDelta\bs'_i,~\text{with } \\
		\|\bDelta\bs'_i\| 
		&\leq \| \bhp_i-\bTheta \bhw_i \| + \| \bhs_i'-\bTheta \bhq_i' \|+ \| \bhS_{i-1} - \bTheta \bhQ_{i-1} \| \| [\bhR]_{(1:i-1,i)}  \| + \|\bTheta \bDelta \bhq_i' \| \\
		&\leq 1.02 u_{fine} n \|\bTheta \|_\Frob   (\|\bhw_i\|+ \|\bhq'_i\|)+  u_{crs}   (0.02\| \bhQ_{i-1}\|+0.01) 1.4 \|\bhw_i \|+4 u_{crs} i^{3/2} \|\bhw_i \| \\
		& \leq 0.02 u_{crs}  \|\bhq'_i\| + 4.05 u_{crs} i^{3/2} \|\bhw_i \| \leq 4.1 u_{crs} i^{3/2} \|\bhw_i \|.
		\end{split}
		\end{equation}  	
		Moreover, it is deduced from~\cref{eq:mainproof1,eq:mainproof3} that 
		\begin{equation} \label{eq:sprime2}
		\begin{split}
		\| \bhs'_i \| &\leq \|\bhp_i\|+\| \bhS_{i-1} \|_\Frob \| [\bhR]_{(1:i-1,i)}  \| + \| \bDelta\bs'_i \| \\
		&\leq 1.25 \|\bhw_i \|  +1.01 i^{1/2} 1.4  \|\bhw_i \| +  4.1 u_{crs} i^{3/2} \|\bhw_i \| \leq 2 i^{3/2}  \|\bhw_i \|.
		\end{split}
		\end{equation}
		
		Consequently, 
		\begin{equation*} 
		\begin{split}
		\|\bhP - \bhS \bhR \|_\Frob 
		&=  \left (\|\bhP_{i-1} - \bhS_{i-1} \bhR_{i-1} \|^2_\Frob + \|\bhp_i - \bhS_{i-1} [\bhR]_{(1:i-1,i)} - \fl(\bhs'_i/\widehat{r}_{i,i}) \widehat{r}_{i,i} \|^2_\Frob \right )^{1/2} \\
		&\leq 4.2  u_{crs}  i^{3/2} (\|\bhW_{i-1}\|_\Frob^2 + \|\bhw_i\|^2)^{1/2}\leq 4.2 u_{crs}   i^{3/2} \|\bhW\|_\Frob, 
		\end{split}
		\end{equation*}
		which results in the first inequality~\cref{eq:maintheorem21} of the theorem.
		
		Furthermore,~\cref{eq:sprime,eq:sprime2} also yield the following results.
		Denote residual matrix $\bhS_{i-1} [\bhR]_{{(1:i-1,1:i-1)}} - \bhP_{i-1}$ by $\bhB_{i-1}$. By the induction hypothesis, $\|\bhB_{i-1} \|_\Frob = \tilde{\Delta}_{i-1} \|\bhP_{i-1} \|_\Frob \leq 4.2 u_{crs} i^{3/2} \| \bhW_{i-1} \|_\Frob$.
		By the standard rounding analysis and the assumption~\cref{eq:Thetanorm}, 
		$$ \sigma_{min}(\bhP_i) \geq \sigma_{min}(\bTheta \bhW_i) - \| \bTheta \bhW_i - \bhP_i \|_\Frob \geq \sigma_{min}(\bTheta \bhW_i) - 0.02 u_{crs} \| \bhW_i\|_\Frob.$$		
		Since $\bTheta$ satisfies the $\varepsilon$-embedding property, we have $\sigma_{min}(\bTheta \bhW_i) \geq \sqrt{1/2}\sigma_{min}(\bhW_i)$. Consequently,
		\begin{equation} \label{eq:sprime3}
		\begin{split}	
		\|\bhs'_i\| &= \|\bhp_i - \bhS_{i-1} [\bhR]_{(1:i-1,i)} + \bDelta\bs'_i \| \\
		& \geq  \|\bhp_i - (\bhP_{i-1} + \bhB_{i-1}) [\bhR]_{{(1:i-1,1:i-1)}}^{-1} [\bhR]_{(1:i-1,i)}\|  - \|\bDelta\bs'_i \| \\
		&\geq \sigma_{min}([\bhP_{i-1} + \bhB_{i-1},\bhp_i]) - \|\bDelta\bs'_i \| \\
		&\geq \sigma_{min}(\bhP_i) - \| \bhB_{i-1}\| - \|\bDelta\bs'_i \| \\
		&\geq \sqrt{1/2}\sigma_{min}(\bhW_i) - 10 u_{crs} i^{3/2} \|\bhW_i\|_\Frob \\
		&\geq 0.68 \sigma_{min}(\bhW_i).	
		\end{split}
		\end{equation}

		In addition, we have 
		\small		
		\begin{equation*}
		\begin{split}	
		&(\bhS_{i-1})^\mathrm{T} \bhs'_i =  (\bhS_{i-1}+ \bDelta \bS_{i-1})^\mathrm{T} \bhs'_i  - (\bDelta \bS_{i-1})^\mathrm{T}\bhs'_i = (\bhS_{i-1}+ \bDelta \bS_{i-1})^\mathrm{T} \bhs'_i  + \bDelta \bt_1 \\ 
		&= (\bhS_{i-1}+\bDelta \bS_{i-1})^\mathrm{T}(\bhp_i - \bhS_{i-1}(\bhS_{i-1}+ \bDelta \bS_{i-1})^\dagger(\bhp_i + \bDelta \bp_i) +  \bDelta\bs'_i) + \bDelta \bt_1 \\
		& = (\bhS_{i-1}+\bDelta \bS_{i-1})^\mathrm{T}(\bhp_i - \bhS_{i-1}(\bhS_{i-1}+ \bDelta \bS_{i-1})^\dagger\bhp_i +  \bDelta\bs'_i) + \bDelta \bt_1 +  \bDelta \bt_2 \\
		& = (\bhS_{i-1}+\bDelta \bS_{i-1})^\mathrm{T}(\bhp_i - (\bhS_{i-1}+ \bDelta \bS_{i-1}) (\bhS_{i-1}+ \bDelta \bS_{i-1})^\dagger\bhp_i )+ \bDelta \bt_1 +  \bDelta \bt_2 + \bDelta \bt_3 + \bDelta \bt_4 \\
		& = \bDelta \bt_1 +  \bDelta \bt_2 + \bDelta \bt_3 + \bDelta \bt_4,
		\end{split}
		\end{equation*}	
		\normalsize
		with 
		\small
		\begin{align*}
		\bDelta \bt_1 
		&= -(\bDelta \bS_{i-1})^\mathrm{T}\bhs'_i, \\
		\bDelta \bt_2 
		&= -(\bhS_{i-1}+\bDelta \bS_{i-1})^\mathrm{T} \bhS_{i-1}(\bhS_{i-1}+ \bDelta \bS_{i-1})^\dagger\bDelta \bp_i,\\
		\bDelta \bt_3 
		&= (\bhS_{i-1}+\bDelta \bS_{i-1})^\mathrm{T} \bDelta \bS_{i-1}(\bhS_{i-1}+ \bDelta \bS_{i-1})^\dagger\bhp_i,\\
		\bDelta \bt_4 
		&= (\bhS_{i-1}+\bDelta \bS_{i-1})^\mathrm{T} \bDelta\bs'_i.
		\end{align*}
		\normalsize
		Since, by the induction hypothesis,
		$$ \|\bhS_{i-1}\| \leq 1.01,~\text{and } \sigma_{min}(\bhS_{i-1}) \geq 0.99, $$
		hence 
		$$ \|\bhS_{i-1} + \bDelta \bS_{i-1} \| \leq \|\bhS_{i-1} \| + \| \bDelta \bS_{i-1} \| \leq  1.021$$
		and 
		$$\|(\bhS_{i-1}+ \bDelta \bS_{i-1})^\dagger\| = 1/(\sigma_{min}(\bhS_{i-1}+ \bDelta \bS_{i-1})) \leq 1/(\sigma_{min}(\bhS_{i-1}) - \|\bDelta \bS_{i-1}\|)  \leq 1.021.  $$
		Consequently, we have
		\small
		\begin{equation} \label{eq:sprime4}
		\begin{split}
		\| \bDelta \bt_1 \| &\leq u_{crs} \|\bhS_{i-1}\| \|\bhs'_i\| \leq 1.01 u_{crs} \|\bhs'_i\|, \\
		\| \bDelta \bt_2 \| &\leq \|\bhS_{i-1}+ \bDelta \bS_{i-1}\| \|\bhS_{i-1} \| \|(\bhS_{i-1}+ \bDelta \bS_{i-1})^\dagger \| \| \bDelta \bp_i \| \leq 0.5 u_{crs} \| \bhp_i  \| \leq  u_{crs} \| \bhw_i \|, \\ 
		\| \bDelta \bt_3 \| &\leq \|\bhS_{i-1}+ \bDelta \bS_{i-1}\| \|\bDelta \bS_{i-1} \| \|(\bhS_{i-1}+ \bDelta \bS_{i-1})^\dagger \|\| \bhp_i \|\leq 0.5 u_{crs} \| \bhp_i \| \leq u_{crs} \| \bhw_i \|, \\
		\| \bDelta \bt_4 \| &\leq \|\bhS_{i-1}+\bDelta \bS_{i-1}\| \|\bDelta \bs_i' \| \leq 4.2 u_{crs} i^{3/2} \| \bhw_i \|.
		\end{split}
		\end{equation}		
		\normalsize
		
		Relations \cref{eq:sprime3,eq:sprime4} imply the following bound 		
		\begin{equation} 
		\begin{split} 
		\|(\bhS_{i-1})^\mathrm{T} \bhs_i \| &\leq  \|(\bhS_{i-1})^\mathrm{T}  (\bhs'_i/\| \bhs'_i \|) \|+0.01 u_{crs} \|\bhS_{i-1}\| \leq  1.03 u_{crs}   + 5 u_{crs} i^{3/2}  \|\bhw_i\|/\|\bhs'_i\|\\
		& \leq 10 u_{crs}  i^{3/2} \frac{\| \bhw_i \|}{\sigma_{min}(\bhW_i)}  \leq  10 u_{crs}  i^{3/2} \cond(\bhW),
		\end{split}
		\end{equation}
		in its turn yielding 	
		\small 
		\begin{equation*} 
		\begin{split}
		\Delta_i 
		&= \|\bI_{i \times i} - \bhS^\mathrm{T} \bhS \|_\Frob =(\|\bI_{(i-1) \times (i-1)} - (\bhS_{i-1})^\mathrm{T} \bhS_{i-1} \|_\Frob^2 + 2\|(\bhS_{i-1})^\mathrm{T} \bhs_i\|^2 +(\|\bhs_i\|^2-1)^2 )^{1/2}\\
		&\leq 20  u_{crs}   i^{2} \cond(\bhW),
		\end{split}
		\end{equation*}
		\normalsize
		which finishes the proof of the theorem for $m=i$.		
	\end{proof}
	
	\begin{proof}[Proof of~\cref{thm:aprioribound}]
		The proof is done by induction on $m$. The statement of the theorem is obvious for $m=1$. Assume that $\bTheta$ is {an} $\varepsilon$-embedding for $\bhQ_m$ for $m=i-1$. As is noted in~\cref{rmk:proofmaintheorem2}, this condition is sufficient for the results in~\cref{thm:maintheorem2} and its proof to hold for $m = i \geq 2$. Therefore we have, 
		\begin{align} \label{eq:theorem2result}
		\|\bhP_i - \bhS_i \bhR_i \|_\Frob 
		&\leq 4.2 u_{crs} i^{3/2} \|\bhW_i \|_\Frob, \\
		\Delta_i = \|\bI_{i \times i} - \bhS_i^\mathrm{T} \bhS_i\|_\Frob 
		&\leq 20  u_{crs} i^2  \cond{(\bhW_i)} \leq 0.02.
		\end{align}	
		And we also have 
		\begin{equation} \label{eq:aprioribound1}
		\|\bhW_i - \bhQ_i \bhR_i \|_\Frob \leq 	{3.3} u_{crs} i^{3/2} \|\bhW_i \|_\Frob.
		\end{equation}
		and 
		\begin{equation} \label{eq:Thetaqs}
		\|\bhS_i - \bTheta \bhQ_i\|_\Frob \leq u_{crs}   (0.02 \|\bhQ_i\| + 0.01), 
		\end{equation}	
		which can be proven similarly to, respectively,{~\cref{eq:qprime} and~\cref{eq:main10}} in the proofs of~\cref{thm:maintheorem1} and~\cref{thm:maintheorem2}. 
		{Now, $$\|\bhP_i - \bTheta \bhW_i \| \leq 1.02 u_{fine} n \| \bTheta \|_\Frob \| \bhW_i \|_\Frob  \leq 1.02 u_{fine} \sqrt{3/2} n^{3/2} i^{1/2} \| \bhW_i \| \leq 0.02 u_{crs}\| \bhW_i \|,$$
			and consequently {$\|\bhS_i \| \leq 1.02$ and $\sigma_{min}(\bhP_i) \geq  \sigma_{min}(\bTheta \bhW_i) - 0.02 u_{crs} \|\bhW_i \|$}, which combined with 
			\cref{eq:theorem2result} and the $\varepsilon$-embedding property of $\bTheta$ results in}
		\begin{equation} \label{eq:singR}
		\sigma_{min}(\bhR_i) \geq \frac{1}{\|\bhS_i\|} (\sigma_{min} (\bhP_i) - \|\bhP_i - \bhS_i \bhR_i \|_\Frob ) \geq 0.7 \sigma_{min}(\bhW_i).
		\end{equation}
		{From~\cref{eq:theorem2result,eq:singR} and the fact that $\|\bhP_i - \bTheta \bhW_i \| \leq 0.02 u_{crs}\| \bhW_i \| $, we deduce that } 
		\begin{equation*} 
		\begin{split}
		\| \bTheta \bhW_i \bhR^{-1}_i - \bhS_i \|_\Frob 
		&\leq \| \bhP_i \bhR^{-1}_i - \bhS_i\|_\Frob + \|(\bhP_i -\bTheta \bhW_i) \bhR^{-1}_i \|_{\Frob} \\
		&\leq (\| \bhP_i - \bhS_i \bhR_i \|_\Frob +  \|\bhP_i -\bTheta \bhW_i\|_\Frob)\|\bhR^{-1}_i\| \\
		& \leq 6.2 u_{crs} i^{2}  \cond(\bhW_i) =: F_1.
		\end{split}
		\end{equation*} 
		This implies that 
		\begin{equation*} 
		\begin{split}
		1  - \Delta_m - F_1 \leq \sigma_{min} (\bTheta \bhW_i \bhR^{-1}_i) \leq \sigma_{max} (\bTheta \bhW_i \bhR^{-1}_i) \leq 1 + \Delta_m+ F_1.
		\end{split}
		\end{equation*} 
		By using the fact that $\bTheta$ is {an} $\varepsilon$-embedding for $\bhW_i$, we arrive to 
		\begin{equation} \label{eq:singWR}
		(1+\varepsilon)^{-1/2} (1  - \Delta_m - F_1) \leq \sigma_{min} (\bhW_i \bhR^{-1}_i) \leq \sigma_{max} (\bhW_i \bhR^{-1}_i) \leq (1-\varepsilon)^{-1/2} (1  + \Delta_m + F_1).
		\end{equation} 			
		We also have from~\cref{eq:aprioribound1,eq:singR},
		\begin{equation*} 
		\|\bhW_i \bhR^{-1}_i - \bhQ_i \|_\Frob \leq 	3.3 u_{crs} i^{3/2} \|\bhW_i \|_\Frob \|\bhR^{-1}_i\|
		\leq 5 u_{crs}  i^{2}  \cond(\bhW_i) =:F_2.
		\end{equation*}
		By combining this relation with~\cref{eq:singWR}, we get
		\begin{equation*} 
		(1+\varepsilon)^{-1/2} (1  - \Delta_m - F_1)-F_2 \leq \sigma_{min} (\bhQ_i) \leq \sigma_{max} (\bhQ_i) \leq (1-\varepsilon)^{-1/2} (1  + \Delta_m + F_1)+F_2.
		\end{equation*} 
		It is deduced that 
		\begin{equation} \label{eq:singQ2}
		(1+\varepsilon)^{-1/2} (1  - F_3) \leq \sigma_{min} (\bhQ_i) \leq \sigma_{max} (\bhQ_i) \leq (1-\varepsilon)^{-1/2} (1  + F_3), 
		\end{equation}
		where $F_3 := \Delta_m + F_1 + \sqrt{5/4}F_2 \leq {32} u_{crs} i^2  \cond{(\bhW_i)}$,
		which in particular implies that $\|\bhQ_i\| \leq 1.6$ and $\sigma_{min} (\bhQ_i) \geq 0.86$. Furthermore, from~\cref{eq:Thetaqs}, we get
		\begin{equation} \label{eq:singthetaQ1}
		1 - F_4 \leq \sigma_{min}(\bhS_i) - 0.1 u_{crs} \leq  \sigma_{min} (\bTheta\bhQ_i) \leq \sigma_{max} (\bTheta\bhQ_i) \leq \sigma_{max}(\bhS_i) + 0.1 u_{crs} \leq  1 + F_4, 
		\end{equation}	
		where $F_4 := \Delta_m +0.1 u_{crs} \leq 20.1 u_{crs} i^2  \cond{(\bhW_i)}$.
		
		Let $\ba \in \mathbb{R}^i$ be some vector. From~\cref{eq:singthetaQ1,eq:singQ2}, we deduce that 
		\begin{align*}
		~~&|  \|\bhQ_i \ba\|^2 - \|\bTheta\bhQ_i \ba\|^2| \\ 
		&\leq \| \ba \|^2 \max \left\{(1-\varepsilon)^{-1}(1+F_3)^2 - (1-F_4)^2    ,(1+F_4)^2 - (1+\varepsilon)^{-1}(1-F_3)^2  \right  \} \\
		&\leq \| \ba \|^2 \max \left\{1.45 \varepsilon+2.9 F_3 +2 F_4, 1.34\varepsilon + 2.1 F_4 + 2.1 F_3\right  \} \\
		& \leq 0.74 \varepsilon^*\| \ba \|^2 \leq \varepsilon^*  \| \bhQ_i \ba \|^2,
		\end{align*}
		where $\varepsilon^* = 2 \varepsilon + 180 u_{crs} i^2 \cond(\bhW_i).$ From this relation and the parallelogram identity, we deduce that $\bTheta$ is a $\varepsilon^*$-embedding for $\bhQ_i$ and finish the proof. 		
	\end{proof}
	\begin{proof}[Proof of~\cref{prop:cert1}]
		First, by standard rounding analysis and assumptions~\cref{eq:Thetanorm,eq:Phinorm}, notice that
		$$ \|\bhV^\bPhi - \bPhi \bV \|_\Frob \leq 1.02 n u_{fine} \|\bPhi\|_\Frob \| \bV \|_\Frob \leq 1.02 \sqrt{\frac{1+\varepsilon}{1-\varepsilon}} u_{fine} n^{3/2} \| \bV^\bPhi \|_\Frob \leq 0.02 u_{crs} \| \bhV^\bPhi \|, $$
		and similarly,
		$$ \|\bhV^\bTheta - \bTheta \bV \|_\Frob \leq 1.02 n u_{fine} \|\bTheta\|_\Frob \| \bV \|_\Frob \leq 1.02 \sqrt{\frac{1+\varepsilon}{1-\varepsilon}} u_{fine} n^{3/2} \| \bV^\bPhi \|_\Frob \leq 0.02 u_{crs} \| \bhV^\bPhi \|. $$
		Furthermore, for any $\bz \in \mathbb{R}^m$, $\|\bz\|=1$, we have	
		\begin{align*}
		\frac{\|\bV^\bTheta \bz \|}{\| \bV^\bPhi \bz \|} 
		&\leq  \frac{\|\bhV^\bTheta \bz \| + \|\bV^\bTheta -\bhV^\bTheta\|}{\| \bhV^\bPhi \bz \|-\|\bV^\bPhi -\bhV^\bPhi\|} \leq \frac{\|\bhV^\bTheta \bz \| + 0.02 u_{crs} \|\bhV^\bPhi\|}{\| \bhV^\bPhi \bz \|-0.02 u_{crs}\|\bhV^\bPhi\|} \\ 
		&\leq (1-0.02 u_{crs}~\cond(\bhV^\bPhi))^{-1} (\frac{\|\bhV^\bTheta \bz \|}{\| \bhV^\bPhi \bz \|} + 0.02 u_{crs}~\cond(\bhV^\bPhi))  ,
		\end{align*}
		and similarly, 
		\begin{align*}
		\frac{\|\bV^\bTheta \bz \|}{\| \bV^\bPhi \bz \|} \geq    (1+0.02 u_{crs}~\cond(\bhV^\bPhi))^{-1} (\frac{\|\bhV^\bTheta \bz \|}{\| \bhV^\bPhi \bz \|} - 0.02 u_{crs}~\cond(\bhV^\bPhi)).
		\end{align*}
		Consequently, 
		$$\left |\frac{\|\bV^\bTheta \bz \|}{\| \bV^\bPhi \bz \|} - \frac{\|\bhV^\bTheta \bz \|}{\| \bhV^\bPhi \bz \|} \right |  \leq 0.02 u_{crs} \cond(\bhV^\bPhi) \left (\frac{\|\bV^\bTheta \bz \|}{\| \bV^\bPhi \bz \|}+1 \right ).  $$
		
		Since 	
		\begin{alignat*}{2}
		\sigma_{min}(\bhV^\bTheta\bhX) &= \min_{\bz \in \mathbb{R}^m, \|\bz\|=1} \frac{\|\bhV^\bTheta \bz \|}{\| \bhV^\bPhi \bz \|}, ~~~ \sigma_{max}(\bhV^\bTheta\bhX) & \shoveleft = \max_{\bz \in \mathbb{R}^m, \|\bz\|=1} \frac{\|\bhV^\bTheta \bz \|}{\| \bhV^\bPhi \bz \|}, \\
		\sigma_{min}(\bV^\bTheta\bX) &= \min_{\bz \in \mathbb{R}^m, \|\bz\|=1} \frac{\|\bV^\bTheta \bz \|}{\| \bV^\bPhi \bz \|}, ~~~ \sigma_{max}(\bV^\bTheta\bX) &= \max_{\bz \in \mathbb{R}^m, \|\bz\|=1} \frac{\|\bV^\bTheta \bz \|}{\| \bV^\bPhi \bz \|}, 
		\end{alignat*}
		we deduce that 
		\begin{align*}
		&|\sigma^2_{min}(\bhV^\bTheta\bhX) - \sigma^2_{min}(\bV^\bTheta\bX)| \leq 0.2 u_{crs} \cond(\bhV^\bPhi)   (\sigma^2_{max}(\bV^\bTheta\bX)+1), \\
		&|\sigma^2_{max}(\bhV^\bTheta\bhX) - \sigma^2_{max}(\bV^\bTheta\bX)| \leq 0.2 u_{crs} \cond(\bhV^\bPhi) (\sigma^2_{max}(\bV^\bTheta\bX)+1).
		\end{align*}
		In its turn, this implies that 
		\begin{align*}
		|\bar{\omega} - \hat{\bar{\omega}} | &\leq 0.2 u_{crs} \cond(\bhV^\bPhi)  (\sigma^2_{max}(\bV^\bTheta\bX)+1) (1+\varepsilon^*) \\ &\leq  0.2 u_{crs} \cond(\bhV^\bPhi)  (\bar{\omega}+1) (1+\varepsilon^*) <u_{crs}~\cond(\bhV^\bPhi). 
		\end{align*}
	\end{proof}

	\begin{proof}[Proof of~\cref{thm:maincorollary1}]
		By~\cref{thm:maintheorem1}, we have the first inequality
		$$  (1+\varepsilon)^{-1/2} (1- \Delta_m - 0.1 u_{crs})  \leq \sigma_{min}(\bhQ) \leq \sigma_{max}(\bhQ) \leq  (1-\varepsilon)^{-1/2} (1+ \Delta_m + 0.1 u_{crs}), $$
		which in particular implies that $\sigma_{min}(\bhQ) \geq 0.7$ and $\sigma_{max}(\bhQ) \leq 1.6$.		
		Furthermore, it also follows from~\cref{thm:maintheorem1} that 	
		$$\|\bhW_m - \bhQ_m \bhR_m\|_\Frob \leq  3.7 u_{crs} m^{3/2} \| \bhW_m\|_\Frob \leq 
		3.7 u_{crs} m^{3/2}  (\| \bb \| + \|[\bhW]_{(1:n,2:m)}\|_\Frob).$$	
		We also have by the standard rounding analysis:
		$$ \|\bhw_i -\bA \bhq_{i-1}\| \leq 1.02 u_{fine} n \|\bA\|_\Frob \|\bhq_{i-1}\| \leq 1.02 u_{fine}  {n^{3/2}} \|\bA \| \|\bhq_{i-1}\| \leq 0.02 u_{crs} m^{-1/2} $$
		These two results imply that 
		$$\bA \bhQ_{m-1} = \bhQ_m  \bhH_{m} + \bDelta \bF_m, \text{ with }$$
		\small
		\begin{equation*}
		\begin{split}
		\| \bDelta \bF_m \|_\Frob 
		&\leq \|[\bhW]_{(1:n,2:m)} - \bA \bhQ_{m-1} \|_\Frob + \|\bhW - \bhQ \bhR\|_\Frob \\
		&\leq 0.02 u_{crs} + 3.7 u_{crs} m^{3/2} (1+\|[\bhW]_{(1:n,2:m)}\|_\Frob)\\
		&\leq 0.02 u_{crs} + 3.7 u_{crs} m^{3/2} (1+\|\bA \| \|\bhQ_{m-1}\|_\Frob + \|[\bhW]_{(1:n,2:m)} - \bA \bhQ_{m-1} \|_\Frob)\\
		&\leq  0.02 u_{crs} + 3.7 u_{crs} m^{3/2} (1+1.6 m^{1/2} +0.02 u_{crs}) 
		\leq 10.5 u_{crs} m^2 
		\end{split}
		\end{equation*}	
		\normalsize
		From this we deduce that 	
		$$\bA \bhQ_{m-1} - \bDelta \bF_m = (\bA - \bDelta \bF_m \bhQ_{m-1}^\dagger) \bhQ_{m-1} = (\bA +\bDelta \bA) \bhQ_{m-1}, $$
		with $ \|\bDelta \bA \|_\Frob = \|\bDelta \bF_m \bhQ_{m-1}^\dagger\|_\Frob \leq \|\bDelta \bF_m\|_\Frob/\sigma_{min}(\bhQ_{m-1}) \leq 10.5 u_{crs} m^2/0.7,$
		which completes the proof.
	\end{proof}

	\begin{proof}[Proof of \Cref{thm:maincorollary2}]
		The proof is done by induction on $m$. The theorem is obvious for $m=1$. Next, the statements are assumed to hold for $m=i-1$ with the objective to prove them for $m=i \geq 2$.	
		
		By induction hypothesis, we have
		$$ \|\bhQ_{i-1}\| \leq 1.6 \text{ and }  \sigma_{min}(\bhQ_{i-1}) \geq 0.7.$$	
		Furthermore, from the relation
		$$ \| \bhW_i - [\bb, \bA \bhQ_{i-1}] \|_\Frob \leq 1.02 u_{fine} {n^{3/2}} \| \bA \| \|\bhQ_{i-1}\|_\Frob \leq u_{crs} \|\bhQ_{i-1}\|_\Frob,$$
		we deduce that  
		$$\cond(\bhW_i) \leq (1+F_1) \cond([\bb, \bA \bhQ_{i-1}]),$$
		where ${F_1 = 2u_{crs} \|\bhQ_{i-1}\|_\Frob/(1-u_{crs} \|\bhQ_{i-1}\|_\Frob) \leq 3.3 u_{crs} m^{1/2} < 0.01}$.		
		Let us next derive a bound for the condition number of matrix $[\bb, \bA \bhQ_{i-1}]$. We have,
		\small
		\begin{align*}
		&\sigma_{min}([\bb, \bA \bhQ_{i-1}]) = \min_{\bv \in \mathbb{R}^i/\{\bnull\}} \frac{\|[\bb, \bA \bhQ_{i-1}]\bv\|}{\|\bv\|} = \min_{\bz_{i-1} \in \mathbb{R}^{i-1}} \frac{\|\bb  - \bA \bhQ_{i-1} \bz_{i-1}\|}{\sqrt{1+\|\bz_{i-1}\|^2}} \\
		& \geq \min_{\bz_{i-1} \in \mathbb{R}^{i-1}} \frac{\|\bb  - \bA \bhQ_{i-1} \bz_{i-1}\|}{1+\|\bz_{i-1}\|} \geq \min_{\bz_{i-1} \in \mathbb{R}^{i-1}} \frac{\|\bb  - \bA \bhQ_{i-1} \bz_{i-1}\|}{1+ \|\bA \bhQ_{i-1} \bz_{i-1}\|/(\sigma_{min}(\bA) \sigma_{min}(\bhQ_{i-1}))} \\
		&\geq \min_{\bz_{i-1} \in \mathbb{R}^{i-1}} \frac{\|\bb  - \bA \bhQ_{i-1} \bz_{i-1}\|}{1+ 1.43 \cond(\bA) \|\bA \bhQ_{i-1} \bz_{i-1}\|} \\
		&\geq  \min_{\bz_{i-1} \in \mathbb{R}^{i-1}} \frac{\|\bb  - \bA \bhQ_{i-1} \bz_{i-1}\|}{1+ 1.43 \cond(\bA)(1 +  \|\bb - \bA \bhQ_{i-1} \bz_{i-1}\|)} 
		\geq  {0.25 {\tau(\bhQ_{i-1})}/\cond(\bA)},	 
		\end{align*}
		\normalsize
		where the last inequality follows due to monotonicity of function $f(x) = x/(a+bx)$ on the interval $[0,+\infty)$, when $a,b>0$, and the fact that $\tau(\bhQ_{i-1}) \leq  \|\bb  - \bA \bhQ_{i-1} \bz_{i-1}\| \leq 1$.
		By combining this fact with
		$ \|[\bb, \bA \bhQ_{i-1}]\| \leq \sqrt{1+ \|\bA \bhQ_{i-1}\|^2} \leq 1.9, $ 
		we conclude that $$\cond(\bhW_i) \leq 8 \cond(\bA)/ \tau(\bhQ_{i-1}).$$ 	
		The theorem is then completed by plugging this result into~\cref{thm:maintheorem2} (with $m=i$).
	\end{proof}
	
	\begin{proof}[Proof of \Cref{thm:GMRESstab}]
		Define $\bDelta \bb = \hat{r}_{1,1} \bhq_1 - \bb$ and notice that $\| \bDelta \bb \| \leq u_{fine}$.  Consider $\bDelta \bA$ from~\cref{thm:maincorollary1}. We have  $\mathrm{range}(\bhQ_{m-1}) = \mathcal{K}_{m-1}(\bA + \bDelta \bA , \bb +\bDelta \bb)$. Furthermore, for any $\bz_{m-1} \in \mathbb{R}^{m-1}$, it holds 
		\begin{equation*}
		\begin{split}
		\| \bhH_m \by_{m-1} - \hat{r}_{1,1} \be_1 \|  
		&\leq \| \bhH_m \bz_{m-1} - \hat{r}_{1,1} \be_1 \| \leq \| \bhQ_m(\bhH_m \bz_{m-1} - \hat{r}_{1,1} \be_1) \|/\sigma_{min}(\bhQ_m) \\
		&= \| (\bA+\bDelta \bA) \bhQ_{m-1}\bz_{m-1} - (\bb +\bDelta \bb) \|/\sigma_{min}(\bhQ_m).
		\end{split}
		\end{equation*}
		On the other hand, 
		\begin{align*}
		\| \bhH_m \by_{m-1} - \hat{r}_{1,1} \be_1 \| 
		&\geq  \| \bhQ_m(\bhH_m \by_{m-1} - \hat{r}_{1,1} \be_1) \|/\|\bhQ_m\| \\
		& = \| (\bA+\bDelta \bA) \bhQ_{m-1}\by_{m-1} - (\bb +\bDelta \bb) \|/\|\bhQ_m\|.
		\end{align*}
		The proposition follows immediately. 
	\end{proof}

\printbibliography[heading=subbibliography]

@article{higham2019new,
	title={A new approach to probabilistic rounding error analysis},
	author={Higham, Nicholas J and Mary, Theo},
	journal={SIAM Journal on Scientific Computing},
	volume={41},
	number={5},
	pages={A2815--A2835},
	year={2019},
	publisher={SIAM}
}

@article{abdelmalek1971round,
	title={Round off error analysis for {G}ram--{S}chmidt method and solution of linear least squares problems},
	author={Abdelmalek, Nabih N},
	journal={BIT Numerical Mathematics},
	volume={11},
	number={4},
	pages={345--367},
	year={1971},
	publisher={Springer}
}

@inproceedings{yamazaki2014mixed,
	title={Mixed-precision orthogonalization scheme and adaptive step size for improving the stability and performance of CA-GMRES on GPUs},
	author={Yamazaki, Ichitaro and Tomov, Stanimire and Dong, Tingxing and Dongarra, Jack},
	booktitle={International Conference on High Performance Computing for Computational Science},
	pages={17--30},
	year={2014},
	organization={Springer}
}

@article{carson2021mixed,
	title={Mixed precision s-step Lanczos and conjugate gradient algorithms},
	author={Carson, Erin and Gergelits, Tom{\'a}{\v{s}} and Yamazaki, Ichitaro},
	journal={Numerical Linear Algebra with Applications},
	pages={e2425},
	year={2021},
	publisher={Wiley Online Library}
}

@article{kim1992efficient,
	title={An efficient parallel algorithm for extreme eigenvalues of sparse nonsymmetric matrices},
	author={Kim, Sun Kyung and Chrortopoulos, AT},
	journal={The International Journal of Supercomputing Applications},
	volume={6},
	number={1},
	pages={98--111},
	year={1992},
	publisher={Sage Publications Sage UK: London, England}
}

@article{barlow2005improved,
	title={Improved {G}ram--{S}chmidt type downdating methods},
	author={Barlow, Jesse L and Smoktunowicz, Alicja and Erbay, Hasan},
	journal={BIT Numerical Mathematics},
	volume={45},
	number={2},
	pages={259--285},
	year={2005},
	publisher={Springer}
}

@article{bjorck1967solving,
	title={Solving linear least squares problems by {G}ram--{S}chmidt orthogonalization},
	author={Bj{\"o}rck, Ake},
	journal={BIT Numerical Mathematics},
	volume={7},
	number={1},
	pages={1--21},
	year={1967},
	publisher={Springer}
}

@article{giraud2005loss,
	title={The loss of orthogonality in the {G}ram--{S}chmidt orthogonalization process},
	author={Giraud, Luc and Langou, Julien and Rozloznik, Miroslav},
	journal={Computers \& Mathematics with Applications},
	volume={50},
	number={7},
	pages={1069--1075},
	year={2005},
	publisher={Elsevier}
}

@article{giraud2005rounding,
	title={Rounding error analysis of the classical {G}ram--{S}chmidt orthogonalization process},
	author={Giraud, Luc and Langou, Julien and Rozlo{\v{z}}n{\'\i}k, Miroslav and van den Eshof, Jasper},
	journal={Numerische Mathematik},
	volume={101},
	number={1},
	pages={87--100},
	year={2005},
	publisher={Springer}
}

@article{leon2013gram,
	title={{G}ram-{S}chmidt orthogonalization: 100 years and more},
	author={Leon, Steven J and Bj{\"o}rck, {\AA}ke and Gander, Walter},
	journal={Numerical Linear Algebra with Applications},
	volume={20},
	number={3},
	pages={492--532},
	year={2013},
	publisher={Wiley Online Library}
}

@article{swirydowicz2018low,
	title={Low synchronization Gram--Schmidt and generalized minimal residual algorithms},
	author={{\'S}wirydowicz, Katarzyna and Langou, Julien and Ananthan, Shreyas and Yang, Ulrike and Thomas, Stephen},
	journal={Numerical Linear Algebra with Applications},
	volume={28},
	number={2},
	pages={e2343},
	year={2021},
	publisher={Wiley Online Library}
}

@inproceedings{malard1994efficiency,
	title={Efficiency and scalability of two parallel QR factorization algorithms},
	author={Malard, J and Paige, CC},
	booktitle={Proceedings of IEEE Scalable High Performance Computing Conference},
	pages={615--622},
	year={1994},
	organization={IEEE}
}

@article{ruhe1983numerical,
	title={Numerical aspects of {G}ram-{S}chmidt orthogonalization of vectors},
	author={Ruhe, Axel},
	journal={Linear algebra and its applications},
	volume={52},
	pages={591--601},
	year={1983},
	publisher={Elsevier}
}

@article{drkovsova1995numerical,
	title={Numerical stability of GMRES},
	author={Drko{\v{s}}ov{\'a}, Jitka and Greenbaum, Anne and Rozlo{\v{z}}n{\'\i}k, M and Strako{\v{s}}, Zdenek},
	journal={BIT Numerical Mathematics},
	volume={35},
	number={3},
	pages={309--330},
	year={1995},
	publisher={Springer}
}

@article{rozloznik1996numerical,
	title={Numerical stability of the GMRES method},
	author={Rozlozn{\'\i}k, Miroslav},
	year={1996},
	publisher={Citeseer}
}

@article{greenbaum1997numerical,
	title={Numerical behaviour of the modified {G}ram-{S}chmidt GMRES implementation},
	author={Greenbaum, Anne and Rozlo{\v{z}}nik, Miroslav and Strako{\v{s}}, Zdenek},
	journal={BIT Numerical Mathematics},
	volume={37},
	number={3},
	pages={706--719},
	year={1997},
	publisher={Springer}
}

@article{paige2006modified,
	title={Modified {G}ram-{S}chmidt (MGS), least squares, and backward stability of MGS-GMRES},
	author={Paige, Christopher C and Rozlozn{\'\i}k, Miroslav and Strakos, Zdenvek},
	journal={SIAM Journal on Matrix Analysis and Applications},
	volume={28},
	number={1},
	pages={264--284},
	year={2006},
	publisher={SIAM}
}

@article{van2004inexact,
	title={Inexact {K}rylov subspace methods for linear systems},
	author={Van Den Eshof, Jasper and Sleijpen, Gerard LG},
	journal={SIAM Journal on Matrix Analysis and Applications},
	volume={26},
	number={1},
	pages={125--153},
	year={2004},
	publisher={SIAM}
}

@article{giraud2007convergence,
	title={Convergence in backward error of relaxed GMRES},
	author={Giraud, Luc and Gratton, Serge and Langou, Julien},
	journal={SIAM Journal on Scientific Computing},
	volume={29},
	number={2},
	pages={710--728},
	year={2007},
	publisher={SIAM}
}

@article{simoncini2007recent,
	title={Recent computational developments in {K}rylov subspace methods for linear systems},
	author={Simoncini, Valeria and Szyld, Daniel B},
	journal={Numerical Linear Algebra with Applications},
	volume={14},
	number={1},
	pages={1--59},
	year={2007},
	publisher={Wiley Online Library}
}

@article{gratton2019exploiting,
	title={Exploiting variable precision in GMRES},
	author={Gratton, Serge and Simon, Ehouarn and Titley-Peloquin, David and Toint, Philippe},
	journal={arXiv preprint arXiv:1907.10550},
	year={2019}
}

@article{ipsen2019probabilistic,
	title={Probabilistic error analysis for inner products},
	author={Ipsen, Ilse CF and Zhou, Hua},
	journal={arXiv preprint arXiv:1906.10465},
	year={2019}
}

@article{higham2020sharper,
	title={Sharper Probabilistic Backward Error Analysis for Basic Linear Algebra Kernels with Random Data},
	author={Higham, Nicholas and Mary, Th{\'e}o},
	year={2020}
}

@article{smoktunowicz2006note,
	title={A note on the error analysis of classical {G}ram--{S}chmidt},
	author={Smoktunowicz, Alicja and Barlow, Jesse L and Langou, Julien},
	journal={Numerische Mathematik},
	volume={105},
	number={2},
	pages={299--313},
	year={2006},
	publisher={Springer}
}

@article{krahmer2011new,
	title={New and improved Johnson--{L}indenstrauss embeddings via the restricted isometry property},
	author={Krahmer, Felix and Ward, Rachel},
	journal={SIAM Journal on Mathematical Analysis},
	volume={43},
	number={3},
	pages={1269--1281},
	year={2011},
	publisher={SIAM}
}

@article{connolly2020stochastic,
	title={Stochastic Rounding and its Probabilistic Backward Error Analysis},
	author={Connolly, Michael P and Higham, Nicholas J and Mary, Th{\'e}o},
	year={2020}
}

@article{pollard17,
	author        = {David Pollard},
	title         = {Lecture notes in Advanced Probability},
	month         = {Spring},
	year          = {2017},
	journal = {Yale university, department of statistics and data science}
	publisher ={Yale university, department of statistics and data science}
}

@book{vershynin2018high,
	title={High-dimensional probability: An introduction with applications in data science},
	author={Vershynin, Roman},
	volume={47},
	year={2018},
	publisher={Cambridge university press}
}

@book{higham2002accuracy,
	title={Accuracy and stability of numerical algorithms},
	author={Higham, Nicholas J},
	year={2002},
	edition ={2},
	publisher={SIAM Publications, Philadelphia, PA, USA}
}

@article{achlioptas2003database,
  title={Database-friendly random projections: Johnson-{L}indenstrauss with binary coins},
  author={Achlioptas, Dimitris},
  journal={Journal of computer and System Sciences},
  volume={66},
  number={4},
  pages={671--687},
  year={2003},
  publisher={Elsevier}
}

@article{ailon2009fast,
  title={Fast dimension reduction using Rademacher series on dual BCH codes},
  author={Ailon, Nir and Liberty, Edo},
  journal={Discrete \& Computational Geometry},
  volume={42},
  number={4},
  pages={615},
  year={2009},
  publisher={Springer}
}

@article{balabanov2019randomized2,
	title={Randomized linear algebra for model reduction. Part {II}: minimal residual methods and dictionary-based approximation},
	author={Balabanov, Oleg and Nouy, Anthony},
	journal={Advances in Computational Mathematics},
	volume={47},
	number={2},
	pages={1--54},
	year={2021},
	publisher={Springer}
}

@article{balabanov2019randomized,
	title={Randomized linear algebra for model reduction. Part {I}: Galerkin methods and error estimation},
	author={Balabanov, Oleg and Nouy, Anthony},
	journal={Advances in Computational Mathematics},
	volume={45},
	number={5-6},
	pages={2969--3019},
	year={2019},
	publisher={Springer}
}

@article{halko2011finding,
  title={Finding structure with randomness: Probabilistic algorithms for constructing approximate matrix decompositions},
  author={Halko, Nathan and Martinsson, Per-Gunnar and Tropp, Joel A},
  journal={SIAM review},
  volume={53},
  number={2},
  pages={217--288},
  year={2011},
  publisher={SIAM}
}

@inproceedings{sarlos2006improved,
	title={Improved approximation algorithms for large matrices via random projections},
	author={Sarlos, Tamas},
	booktitle={Foundations of Computer Science, 2006. FOCS'06. 47th Annual IEEE Symposium on},
	pages={143--152},
	year={2006},
	organization={IEEE}
}

@article{woodruff2014sketching,
	title={Sketching as a tool for numerical linear algebra},
	author={Woodruff, David P and others},
	journal={Foundations and Trends{\textregistered} in Theoretical Computer Science},
	volume={10},
	number={1--2},
	pages={1--157},
	year={2014},
	publisher={Now Publishers, Inc.}
}
\end{refsection}
\end{document}